\documentclass[12pt,reqno]{amsart}
\usepackage[numbers,sort&compress]{natbib}
\usepackage{amsfonts}
\usepackage{amssymb,color}
\usepackage{amsmath}

\usepackage[colorlinks=true]{hyperref}
\hypersetup{urlcolor=blue, citecolor=blue, linkcolor=blue}

\usepackage[margin=3cm, a4paper]{geometry}

\makeatletter
\renewcommand{\section}{\@startsection{section}{1}{\z@}%
  {1.4\linespacing\@plus.3\linespacing\@minus.2\linespacing}%
  {.7\linespacing\@plus.1\linespacing\@minus.1\linespacing}%
  {\normalfont\scshape\centering}}
\makeatother

\newtheorem{thm}{Theorem}[section]
\newtheorem{cor}[thm]{Corollary}
\newtheorem{lem}[thm]{Lemma}
\newtheorem{prop}[thm]{Proposition}
\theoremstyle{definition}

\theoremstyle{plain}

\theoremstyle{remark}
\newtheorem{remark}[thm]{Remark}

\newcommand{\EQ}[1]{\begin{align*}\begin{split} #1 \end{split}\end{align*}}
\newcommand{\EQn}[1]{\begin{align}\begin{split} #1 \end{split}\end{align}}
\newcommand{\EQnn}[1]{\begin{align} #1 \end{align}}

\def\norm#1{\left\|#1\right\|}
\def\normo#1{\|#1\|}
\def\normb#1{\big\|#1\big\|}

\def\abs#1{\left|#1\right|}

\def\brk#1{\left(#1\right)}

\def\brkb#1{\big(#1\big)}

\def\mbrk#1{\left[#1\right]}

\def\mbrkb#1{\big[#1\big]}
\def\mbrkbb#1{\bigg[#1\bigg]}

\def\fbrk#1{\left\lbrace#1\right\rbrace}

\def\jb#1{\langle#1\rangle}

\def\wt#1{\widetilde{#1}}
\def\wh#1{\widehat{#1}}
\def\wb#1{\overline{#1}}

\def\lsm{\lesssim}
\def\gsm{\gtrsim}

\newcommand{\R}{{\mathbb R}}
\newcommand{\C}{{\mathbb C}}

\newcommand{\F}{{\mathcal{F}}}

\newcommand{\B}{{\mathcal{B}}}

\newcommand{\K}{{\mathcal{K}}}
\newcommand{\M}{{\mathcal{M}}}

\newcommand{\U}{{\mathcal{U}}}
\newcommand{\V}{{\mathcal{V}}}

\newcommand{\dd}{{\mathrm{d}}}

\def\Om{\Omega}

\def\ph{\varphi}

\def\ga{\gamma}

\newcommand{\I}{\infty}

\numberwithin{equation}{section}

\begin{document}
\title[Critical long-range NLS]{Modified wave operators for the critical long-range nonlinear Schr\"odinger equation}
\subjclass[2020]{35Q55, 35B40, 35B30}
\keywords{Nonlinear Schr\"odinger equation, long-range scattering, modified wave operator, large final state, linearized equation}
\author[Jia Shen]{Jia Shen}
\address{(J. Shen) School of Mathematical Sciences and LPMC\\
	Nankai University\\ Tianjin 300071, China}
\email{shenjia@nankai.edu.cn}
\thanks{}

\author{Yifei Wu}
\address{(Y.Wu) School of Mathematical Sciences\\
	Nanjing Normal University\\
	Nanjing 210046, China}
\email{yerfmath@gmail.com}
\date{}

\begin{abstract}
In this paper, we construct modified wave operators for the critical long-range nonlinear Schr\"odinger equation with final data $\jb{x}^su_+\in L_x^2,\ d/2<s\le2$, for dimensions $d=1,2,3$, and for both defocusing and focusing cases. We impose no smallness, nonvanishing, or spatial derivative assumption. Most previous modified wave operator results require smallness on $\norm{\wh u_+}_{L^\I}$, analytic final data or smooth amplitude-phase structure, and the large data modified wave operator problem with finite weighted condition is only solved for 1D defocusing equation.

We introduce a new perturbation method which successfully applies to both defocusing and focusing signs when $d\le3$. The main ingredient is an improved growth estimate for the linearized flow, obtained by a low--high frequency decomposition. At low frequencies, spatial localization and freezing the coefficients in space lead to a Fuchsian system. The nilpotent structure at the singular point and a Hamiltonian energy estimate give logarithmic growth for this system. At high frequencies, a modified energy for suitable symmetry operator cancels the leading contribution to energy growth.  
\end{abstract}
\maketitle
\enlargethispage{2pt}

\section{Introduction}\label{sec:intro}

We consider the nonlinear Schr\"odinger equation (NLS)
\EQn{\label{eq:nls}
i\partial_tu+\frac12\Delta u=\lambda|u|^pu,
\qquad (t,x)\in\R\times\R^d.
}
Throughout the paper, unless otherwise stated, we assume
\EQn{\label{eq:parameters}
d\in\{1,2,3\},\quad p:=\frac2d,\quad \lambda\in\{-1,1\}.
}
The final-state problem is long range, and the asymptotic profile requires a logarithmic phase correction. We define the weighted space
\EQn{\label{lr:weighted-space}
\F H^s(\R^d):=\{f\in L^2(\R^d):\jb{ x}^sf\in L^2(\R^d)\},
\qquad \jb{ x}:=(1+|x|^2)^{1/2}.
}

The failure of linear scattering in the long-range regime goes back to Strauss~\cite{Strauss1974} and Barab~\cite{Barab1984}. The construction of a modified wave operator is a final-state problem: one prescribes $u_+$ and seeks a nonlinear solution with the prescribed modified asymptotics. The corresponding completeness problem is that given suitable initial data, construct a nonlinear solution such that the modified scattering holds.

For the critical final-state problem, Ozawa~\cite{Ozawa1991} first constructed modified wave operators for the 1D cubic NLS. Ginibre and Ozawa~\cite{GinibreOzawa1993} extended the construction to $d=2,3$. Carles~\cite{Carles2001} gave a geometric-optics construction in one dimension, and Hayashi and Naumkin~\cite{HayashiNaumkin2006} subsequently weakened the weighted regularity assumptions. All these results require suitable regularity and decay of $u_+$, together with the smallness condition $
\normo{\wh{u_+}}_\infty\ll1$.

For the initial-value problem, Hayashi and Naumkin~\cite{HayashiNaumkin1998} proved modified scattering for small localized initial data in $d=1,2,3$. Other methods have been exploited for 1D cubic case: Lindblad and Soffer~\cite{LindbladSoffer2006} used an asymptotic ODE in physical space to prove completeness for small initial data. Kato and Pusateri~\cite{KatoPusateri2011} gave a proof using the space time resonance method, and Ifrim and Tataru~\cite{IfrimTataru2014,IfrimTataru2024WavePackets} developed the method of testing by wave packets.

Large-data asymptotics have also been studied for the initial data problem. For the integrable 1D defocusing cubic NLS, Deift and Zhou~\cite{DeiftZhou2003} obtained long-time asymptotics for arbitrary initial data in a weighted Sobolev space by the Riemann--Hilbert steepest descent method. Cazenave and Naumkin~\cite{CazenaveNaumkin2018} constructed modified-scattering solutions in arbitrary dimension from sufficiently oscillatory initial data satisfying weighted regularity and a quantitative nonvanishing condition. Georgiev and Ozawa~\cite{GeorgievOzawa2025} obtained a conditional completeness result for the defocusing critical equation in dimensions one and two, assuming $u_0\in \F H^1$ and $\normo{u(t)}_\infty\lsm t^{-d/2}$ for $t\ge1$.

Earlier constructions for large final data impose stronger conditions on $u_+$. As part of their Gevrey-space theory, Ginibre and Velo~\cite{GinibreVelo2001} constructed modified wave operators for both signs of the one-dimensional cubic equation without a size restriction. A sufficient condition on the final datum is
\EQ{
e^{\rho|x|}\jb{ x}^3u_+\in L^2(\R)
\qquad\text{for some }\rho>0,
}
which gives analytic regularity of $\wh{u_+}$. For 1D cubic equation, Lindblad and Soffer~\cite{LindbladSoffer2006} constructed solutions by higher-order asymptotic expansions. Their final data admit a representation $\wh{u_+}=a e^{ib}$ with real-valued, smooth, rapidly decaying amplitude and phase. In particular, $a,b\in\mathcal S(\R;\R)$ give a sufficient class, with no amplitude restriction in the defocusing case.

More recently, Kawamoto and Mizutani~\cite{KawamotoMizutani2025} constructed modified wave operators for the 1D defocusing cubic equation for arbitrary $
u_+\in \F H^{1+\varepsilon}(\R),\ \varepsilon>0$. 
Their focusing result additionally requires $|\lambda|\normo{\wh{u_+}}_\infty^2<1$. Therefore, the large data modified wave operators problem on the finite weighted space is left unsolved for 1D focusing case, and for higher dimensional cases.

For subcritical powers $0<p<2/d$, Hayashi, Kaikina, and Naumkin~\cite{HayashiKaikinaNaumkin1999} studied the 1D Cauchy problem for small initial data of the form $u(1,x)=e^{ix^2/2}v_0(x)$, where $v_0$ belongs to polynomially weighted analytic spaces and satisfies a quantitative polynomial lower bound. For the defocusing NLS, they obtained an explicit ODE-type asymptotic profile when $1<p<2$. For $d\ge2$, Han~\cite{Han2024} constructed global solutions with linear pointwise decay from a class of arbitrarily large initial data with a sufficiently strong quadratic phase, weighted regularity, and a quantitative nonvanishing condition. He obtained phase-modified asymptotics in the range $
\frac4{d+\sqrt{d^2+8d}}<p<\frac2d$. 

Very recently, there are some new developments towards the modified wave operator problem for long range NLS. The authors~\cite{ShenWu2026Subcritical} constructed modified wave operators for subcritical long range NLS throughout $0<p<2/d$ in all dimensions $d\ge1$, for final data with small, analytic, nonvanishing, and logarithmic boundedness assumptions. Particularly, when $0<p\le1/d$, the leading ODE phase correction is not sufficient, and higher-order amplitude and phase corrections are identified from a Fuchsian system, and a nonlinear final-state normal form resolves the remaining long-range coupling. On the other direction, Kawamoto and Mizutani~\cite{KawamotoMizutani2026Subcritical} constructed large-data modified wave operators for the 1D defocusing equation when $4/\sqrt7<p<2$. Their result applies to every prescribed $u_+\in \F H^3(\R)$, and also allows lower weights depending on $p$.

In the present critical long range problem, we construct modified wave operators for arbitrary $u_+\in \F H^s$ with $d/2<s\le2$ and both defocusing and focusing cases. For a prescribed final datum $u_+$, set $W:=\wh{u_+}$ and $B:=\normo{W}_\infty$. For $t\ge2$, define
\EQn{\label{eq:profile}
V(t,y)&:=W(y)\exp\brkb{-i\lambda|W(y)|^p\log t},\\
A_W(t,x)&:=\frac1{(it)^{d/2}}e^{i|x|^2/(2t)}
V\brk{t,\frac xt}.
}

\begin{thm}[Modified wave operators for large final states]
\label{thm:main}
Let $u_+\in \F H^s(\R^d)$ with $d/2<s\le2$. Then there exists a unique global $L^2$ solution $u$ of~\eqref{eq:nls} such that, for all sufficiently large $t$,
\EQn{\label{eq:main-rate}
\normo{u(t)-A_W(t)}_{L_x^2}
+\normo{u-A_W}_{L_{\tau,x}^{\frac{2(d+2)}d}([t,\infty)\times\R^d)} \lsm_{u_+,s}t^{-s/2}(\log t)^2.
}
This defines the modified wave operator
\EQ{
\Om_+:\F H^s(\R^d)\longrightarrow L^2(\R^d),
\qquad \Om_+(u_+):=u(0).
}
\end{thm}

\begin{remark}[Why smallness is required previously]\label{rem:terminal-integral}
The construction of a modified wave operator requires solving a nonlinear final data problem around the prescribed asymptotic, which leads to integrating the forced term of some linearized equation. Heuristically, the desired forcing has decay $r^{-1-a}$ with $a>0$, and the linearized equation has the growth bound $(r/t)^{\kappa_B}$ backward from $r$ to $t$. Then, the integral
	\EQ{
		\int_t^\infty\brk{\frac rt}^{\kappa_B}r^{-1-a}\,\dd r,
	}
should be bounded. The direct coefficient bound for the linearization around $A_W$ gives $\kappa_B:=(1+p)B^p$, where $B=\norm{\wh u_+}_{L^\I}$. In the usual small-data argument based on this bound, one imposes $(1+p)B^p<a$. Increasing the starting time does not change this restriction. For more details, we refer the reader to Section \ref{sec:main-difficulty} below. For the 1D defocusing cubic equation, Kawamoto and Mizutani~\cite{KawamotoMizutani2025} overcome this obstruction by using energy structure based on pseudo-conformal transformation. 
\end{remark}

\begin{remark}[Main idea of the method]
The main point is to improve the time growth of the linearized flow for arbitrary amplitudes. A low--high frequency decomposition separates two mechanisms:
\begin{itemize}
	\item (Low frequency part)
	The nilpotent and Hamiltonian structures of the Fuchsian system obtained after freezing the spatial coefficients.
	\item (High frequency part)
	Modified energy estimate at high frequencies.
\end{itemize}
Together they reduce the growth exponent below the decay exponent required for the final-state integral. The nonlinear remainder then yields small factors at large times. This issue is discussed in more details in Section \ref{sec:proof-ideas} below. 
\end{remark}

\begin{remark}
The combination of nilpotent resonant dynamics and modified energy estimates developed here suggests a general perturbative framework for constructing solutions with large prescribed asymptotic profiles. The method has the potential to solve the large final data problems for the following models for instance:
\begin{itemize}
	\item Derivative nonlinear Schrödinger equation;
	\item Critical long-range nonlinear Klein–Gordon equation;
	\item Modified Korteweg–de Vries equation. (It would require additional analysis of the third harmonic generated by cubic interactions, as well as a treatment of the self-similar region and its connection with the oscillatory asymptotics.)
\end{itemize} 
\end{remark}

\begin{remark}[The role of heat regularization]\label{rem:heat-regularization}
The above mentioned method can be carried out directly for $u_+\in \F H^2$. To further lowering the weight, we use the heat regularization to relax the assumption to $u_+\in \F H^s$ with $d/2<s<2$. If $u_+\in \F H^2$, the proof can be further simplified.
\end{remark}


\subsection{The main difficulty}\label{sec:main-difficulty}

Write $W=\wh{u_+}$ and $B=\normo{W}_\infty$. The linearized equation for the correction $h=u-A_W$ is
\EQn{\label{eq:linearized-physical}
\brk{i\partial_t+\frac12\Delta}h
=\lambda\mbrk{
\brk{1+\frac p2}|A_W|^ph
+\frac p2|A_W|^{p-2}A_W^2\wb h
}+f.
}
The forcing $f$ contains the nonlinear remainder and the profile error. To isolate the difficulty in the linear estimate, consider this equation where the forcing has decay $\normo{f(r)}_2\lsm r^{-1-a}$ with some $a>0$. The linearized potential has size $t^{-1}$:
\EQ{
|A_W(t,x)|^p=t^{-1}\abs{W\brk{\frac xt}}^p,
\qquad
\normo{|A_W(t)|^p}_\infty=\frac{B^p}{t}.
}
A direct Gronwall estimate therefore gives for $r\ge t$:
\EQn{\label{eq:intro-gronwall}
\normo{h(t)}_2
&\le\brk{\frac rt}^{\kappa_B}\normo{h(r)}_2
+\int_t^r\brk{\frac{\tau}{t}}^{\kappa_B}\normo{f(\tau)}_2\,\dd \tau,
}
where $\kappa_B=(1+p)B^p$. Thus the growth exponent depends on the size $\normo{\wh{u_+}}_\infty$ of the final datum. To construct a solution, one must estimate the forcing integral up to infinity. However, its upper bound is
\EQ{
\int_t^\infty\brk{\frac rt}^{\kappa_B}r^{-1-a}\,\dd r.
}
This direct estimate requires $(1+p)B^p<a$, while increasing the starting time does not remove this restriction.

\subsection{Our approach}\label{sec:proof-ideas}

Our task is to replace $\kappa_B$ in~\eqref{eq:intro-gronwall} by an arbitrarily small exponent. By Duhamel's formula, it suffices to estimate solutions of~\eqref{eq:linearized-physical} with $f=0$. We achieve this using the following strategy.

\subsubsection{The ray equation}

For such a solution $h$, we introduce $z$ by
\EQ{
h(t,x)=\frac{e^{i|x|^2/(2t)}}{(it)^{d/2}}z\brk{t,\frac xt}.
}
The resulting linearized ray equation is
\EQn{\label{eq:intro-ray-equation}
i\partial_tz+\frac1{2t^2}\Delta_yz
&=\frac\lambda t\mbrk{
\brk{1+\frac p2}|W|^pz+\frac p2|W|^{p-2}V(t)^2\wb z
},\\
V(t)&=We^{-i\lambda|W|^p\log t}.
}
We use this change of variables for the following reasons.
\begin{itemize}
\item It removes the quadratic oscillation and spatial dilation in $A_W$. The profile $W(x/t)$ becomes the fixed spatial profile $W(y)$, which is suited to localization and freezing of coefficients.
\item The dispersive term has coefficient $1/(2t^2)$, while the linearized potential has coefficient $1/t$. These distinct time scales allow us to identify the Fuchsian system after freezing the coefficients and to isolate the dominant skew-adjoint part at high frequencies.
\end{itemize}
Thus the problem reduces to proving that, for every $\eta>0$, solutions of~\eqref{eq:intro-ray-equation} satisfy
\EQn{\label{eq:intro-linear}
\normo{z(t)}_2\lsm_{W,\eta}
\brk{\frac rt}^{\eta}\normo{z(r)}_2,
\qquad r\ge t\ge T_{W,\eta}.
}
Taking real and imaginary parts, set
\EQ{
X(t,y):=\begin{pmatrix}
\operatorname{Re}z(t,y)\\
\operatorname{Im}z(t,y)
\end{pmatrix},
\qquad
J_0:=\begin{pmatrix}0&-1\\1&0\end{pmatrix}.
}
The linearized ray equation~\eqref{eq:intro-ray-equation} becomes the real system
\EQn{\label{eq:intro-real-system}
\partial_tX
&=\frac1{2t^2}J_0\Delta_yX+\frac1t\brk{\V+\B(t)}X,\\
\V&=-\lambda\brk{1+\frac p2}|W|^pJ_0,\quad
\B(t)=\frac{\lambda p}{2}|W|^{p-2}
\begin{pmatrix}
\operatorname{Im}(V^2)&-\operatorname{Re}(V^2)\\
-\operatorname{Re}(V^2)&-\operatorname{Im}(V^2)
\end{pmatrix}.
}
This real formulation will be used both to identify the Fuchsian system after freezing the coefficients and to construct the modified energy at high frequencies.

On each time interval, we use a smooth Fourier decomposition $z(t)=P_{\le N}z(t)+P_{>N}z(t)$, where $P_{\le N}$ is the low-frequency cutoff and $P_{>N}:=I-P_{\le N}$. 

\subsubsection{Low frequencies}

To estimate $P_{\le N}X$, we apply spatial localization to its equation obtained from the real system~\eqref{eq:intro-real-system}, and freeze $W$ at a constant $w_0$ with $|w_0|\le B$. For the homogeneous part of this frozen system, a phase rotation and the Fourier transform in $y$, followed by $r:=|\xi|^2/(2t)$ for $\xi\ne0$, give the Fuchsian system for $\vec u=(u_1,u_2)^T$:
\EQ{
\partial_r\vec u
=\brk{\frac{A_g}{r}+J_0}\vec u,\quad 
A_g:=\begin{pmatrix}0&0\\g&0\end{pmatrix},
\quad g:=\lambda p|w_0|^p.
}
For small $r$, the nilpotent structure gives only logarithmic growth, and variation of constants controls the bounded term $J_0$ on a fixed small-$r$ interval. The improvement in time growth due to the nilpotent structure was already observed by Banica and Vega~\cite{BanicaVega2012,BanicaVega2013} for linearized equations with spatially constant coefficients. 

For large $r$, the Hamiltonian structure gives the positive energy
\EQ{
E(r)&:=(1+g/r)|u_1(r)|^2+|u_2(r)|^2
\sim\normo{\vec u(r)}^2,\\
E'(r)&=-\frac g{r^2}|u_1(r)|^2,\qquad r\ge2(1+|g|).
}
Its growth rate is integrable for either sign of $g$. Combining the two regions yields
\EQ{
\normo{\vec u(r_1)}
\lsm_{p,B}\brk{1+\abs{\log\frac{r_1}{r_0}}}\normo{\vec u(r_0)},
}
which is sufficient for~\eqref{eq:intro-linear} for every $\eta>0$ in the low frequency part. However, if we do not use the above strategy, Gronwall's inequality applied directly to Duhamel's formula only gives the power bound $\brk{\frac{\max\{r_0,r_1\}}{\min\{r_0,r_1\}}}^{|g|} $.

\subsubsection{High frequencies}

Let $w:=P_{>N}X$. We mollify the coefficients in $y$ and keep the notation $\V,\B$ for the resulting matrices in this calculation. Then
\EQ{
\partial_tw
&=\brk{\frac1{2t^2}J_0\Delta_y+\frac1t\V}w
+\frac1t\B(t)w+f_h,\\
\V^*&=-\V,\qquad
\B(t)^*=\B(t),\qquad
J_0\B(t)=-\B(t)J_0,
}
where $f_h$ contains the coefficient approximation and Fourier cutoff errors. Without the energy correction, the symmetric term gives
\EQ{
\abs{\partial_t\normo{w}_2^2}
\lsm_p\frac{B^p}{t}\normo{w}_2^2
       +\normo{f_h}_2\normo{w}_2.
}

We define 
\EQ{
\K(t):=tJ_0\B(t)(-\Delta_y)^{-1}P_{>N/2},
}
and the leading symmetric term is absorbed into the derivative of the modified energy
\EQ{
E(t):=\normo{w(t)}_2^2+2\jb{\K(t)w(t),w(t)}
\sim\normo{w(t)}_2^2.
}
On each finite time interval, at sufficiently high frequencies, we obtain for any prescribed $\eta>0$
\EQ{
\abs{E'(t)}
\lsm_{p,B}\frac{\eta}{t}E(t)
       +\normo{f_h(t)}_2E(t)^{1/2}.
}
The improvement is in the coefficient of $t^{-1}$: the fixed amplitude-dependent coefficient $B^p$ is replaced by an arbitrarily small coefficient $\eta$, while $E$ remains equivalent to $\normo{w}_2^2$. The coefficient approximation and cutoff errors in $f_h$ are controlled at large times. This gives the high-frequency estimate needed for~\eqref{eq:intro-linear}.

\section{Preliminaries}\label{sec:preliminaries}

\subsection{Notation}\label{sec:notation}

We write $A\lsm B$ or $B\gsm A$ if $A\le CB$ for some positive constant $C$. The constant $C$ may change from line to line. Parameter dependence is indicated by subscripts. We write $\normo{\cdot}_r$ for the norm on $L^r(\R^d)$ and $D^2f$ for the full spatial Hessian, equipped with its Euclidean norm. We use the usual inhomogeneous Sobolev spaces $H^s(\R^d)$.  

We use the Fourier transform
\EQn{\label{lr:fourier-convention}
 \wh f(\xi):=\mathcal Ff(\xi):=(2\pi)^{-d/2}\int_{\R^d}e^{-ix\cdot\xi}f(x)\,\dd x.
}
For $t>0$, the change of variables and its inverse are
\EQn{\label{eq:ray-map}
 (\mathcal P_tf)(x):=(it)^{-d/2} e^{i|x|^2/(2t)} f\brk{\frac xt};\quad
 (\mathcal P_t^{-1}g)(y)=(it)^{d/2}e^{-it|y|^2/2}g(ty).
}

Differentials of maps from $\C$ to $\C$ are taken over $\R$. We write $F(z):=|z|^pz$. Its real differential $DF(z)$ is the real-linear map given by
\EQ{
DF(z)[h]=\brk{1+\frac p2}|z|^ph
 +\frac p2|z|^{p-2}z^2\wb h,
 \qquad z,h\in\C.
}
We define the nonlinear remainder by
\EQ{
\mathcal N_z(h):=F(z+h)-F(z)-DF(z)[h].
}
Throughout the paper, scalar and matrix coefficients with removable singularities at zero, such as $|z|^{p-2}z^2$, are understood by continuous extension.

For a compact time interval $I$, let
\EQn{\label{eq:nl-S}
\normo{v}_{S(I)}
 :=\normo{v}_{L^\infty(I;L^2)}
  +\normo{v}_{L^{\frac{2(d+2)}d}(I\times\mathbb R^d)}.
}

Let $\chi\in C_c^\infty(\R^d)$ be real and radial, with $0\le\chi\le1$, $\chi=1$ on $\{|\zeta|\le1\}$, and $\chi=0$ on $\{|\zeta|\ge2\}$. For $N>0$, set
\EQn{\label{lin:P-Q}
P_{\le N}:=\chi(D_y/N),\qquad P_{>N}:=I-P_{\le N},\qquad D_y:=-i\nabla_y.
}
For $0<\gamma\le1$ and a bounded $\gamma$-H\"older continuous matrix multiplier $b$, we use the standard commutator estimate
\EQn{\label{lin:P-commutator}
\normo{[P_{\le N},b]f}_2+\normo{[P_{>N},b]f}_2
 \lsm_\gamma N^{-\gamma}[b]_{C^{0,\gamma}}\normo{f}_2.
}
Fix $\ph\in C_c^\infty(\R^d;\R)$, 
 $\ph\ge0$,  $\int_{\R^d}\ph(x)\,\dd x=1$. 
For $\theta>1$, define
\EQn{\label{lin:mollifier}
\ph_\theta(x):=e^{d\theta/4}\ph(e^{\theta/4}x),\qquad
 b_\theta:=\ph_\theta*b.
}
Then $\normo{b_\theta}_\infty\le\normo{b}_\infty$, and the standard mollification estimates give
\EQn{\label{lin:mollification-estimates}
\normo{b-b_\theta}_\infty\lsm_\gamma e^{-\gamma\theta/4}[b]_{C^{0,\gamma}},\qquad
 \normo{\nabla^jb_\theta}_\infty\lsm_j e^{j\theta/4}\normo{b}_\infty,
 \qquad j=1,2.
}

\subsection{Useful lemmas}

We use the following weighted interpolation estimate to control the gradient term that occurs when differentiating the phase map twice.

\begin{lem}[Composition estimate]\label{lem:weighted-second}
Let $f\in H^2(\R^d)$. Define $|f|^{p-1}|\nabla f|^2$ to be zero on $\{f=0\}$. Then
\EQn{\label{eq:weighted-second}
\normb{|f|^{p-1}|\nabla f|^2}_2
 \lsm_p \normo{f}_\infty^p\normo{D^2f}_2.
}
\end{lem}

\begin{proof}
For $\varepsilon>0$, set $f_\varepsilon:=(|f|^2+\varepsilon^2)^{1/2}$. We first prove
\EQn{\label{eq:weighted-regularized}
\normb{f_\varepsilon^{p-1}|\nabla f|^2}_2
 \lsm_p (\normo{f}_\infty^2+\varepsilon^2)^{p/2}\normo{D^2f}_2.
}
For $f\in C_c^\infty(\R^d)$, set
\EQ{
 I_\varepsilon:=\int_{\R^d}f_\varepsilon^{2p-2}|\nabla f|^4.
}
Integrating the divergence of $f_\varepsilon^{2p-2}|\nabla f|^2\operatorname{Re}(\wb f\nabla f)$ gives
\EQ{
&I_\varepsilon+2(p-1)\int f_\varepsilon^{2p-4}
       |\operatorname{Re}(\wb f\nabla f)|^2|\nabla f|^2\\
 &\quad=-\int f_\varepsilon^{2p-2}|\nabla f|^2
                  \operatorname{Re}(\wb f\Delta f) -\int f_\varepsilon^{2p-2}
       \operatorname{Re}(\wb f\nabla f)\cdot\nabla(|\nabla f|^2).
}
By $|f|\le f_\varepsilon$,
\EQ{
|\operatorname{Re}(\wb f\nabla f)|^2
\le |f|^2|\nabla f|^2
\le f_\varepsilon^2|\nabla f|^2.
}
Since $p=2/d\ge2/3$, this gives
\EQ{
I_\varepsilon+2(p-1)\int f_\varepsilon^{2p-4}
|\operatorname{Re}(\wb f\nabla f)|^2|\nabla f|^2
\ge\min\{1,2p-1\}I_\varepsilon
\ge\tfrac13 I_\varepsilon.
}
Using $|\nabla(|\nabla f|^2)|\le2|\nabla f||D^2f|$, we obtain
\EQn{\label{eq:weighted-ibp}
I_\varepsilon
 \lsm \int f_\varepsilon^{2p-1}|\nabla f|^2|D^2f| \lsm I_\varepsilon^{1/2}
       \brk{\int f_\varepsilon^{2p}|D^2f|^2}^{1/2} \lsm I_\varepsilon^{1/2}
       (\normo{f}_\infty^2+\varepsilon^2)^{p/2}\normo{D^2f}_2.
}
Dividing by $I_\varepsilon^{1/2}$ when it is nonzero proves~\eqref{eq:weighted-regularized}. For fixed $\varepsilon>0$, density and the embeddings $H^2\hookrightarrow L^\infty\cap W^{1,4}$ extend the estimate to $f\in H^2$, and Fatou's lemma on $\{f\ne0\}$ as $\varepsilon\downarrow0$ gives~\eqref{eq:weighted-second}.
\end{proof}

For $\ell\ge0$, set $G_\ell(z):=z e^{-i\lambda\ell|z|^p}$. This is a real $C^1$ map on $\C$. For $z,h\in\C$, the real derivative is
\EQn{\label{lr:phase-first-derivative}
DG_\ell(z)[h]
 =e^{-i\lambda\ell|z|^p}
    \brk{h-i\lambda\ell p|z|^{p-2}z
                         \operatorname{Re}(\wb zh)}.
}
We also have
\EQ{
\normo{DG_\ell(z)}\le1+p\ell|z|^p,
 \qquad
 \partial_\ell G_\ell(z)=-i\lambda|z|^pG_\ell(z).
}

\begin{lem}\label{lem:phase-hessian}
Let $f\in H^2(\R^d)$ and $\ell\ge0$. Then $G_\ell(f)\in H^2$ and
\EQn{\label{eq:phase-hessian}
\normo{D^2G_\ell(f)}_2
 \lsm_p (1+\ell\normo{f}_\infty^p)^2\normo{D^2f}_2.
}
\end{lem}

\begin{proof}
For $\varepsilon>0$, set
\EQ{
f_\varepsilon:=(|f|^2+\varepsilon^2)^{1/2},\qquad
G_{\ell,\varepsilon}(f):=f e^{-i\lambda\ell f_\varepsilon^p}.
}
Differentiating twice gives
\EQ{
|D^2G_{\ell,\varepsilon}(f)|
 \lsm_p (1+\ell f_\varepsilon^p)|D^2f| +\ell(1+\ell f_\varepsilon^p)
       f_\varepsilon^{p-1}|\nabla f|^2.
}
Taking $L^2$ norms and applying~\eqref{eq:weighted-regularized} from the proof of Lemma~\ref{lem:weighted-second}, we obtain
\EQ{
\normo{D^2G_{\ell,\varepsilon}(f)}_2
\lsm_p \mbrk{1+\ell(\normo{f}_\infty^2+\varepsilon^2)^{p/2}}^2
       \normo{D^2f}_2.
}
For $p<1$, we have
\EQ{
f_\varepsilon^{p-1}|\nabla f|^2
\le |f|^{p-1}|\nabla f|^2\in L^2,
}
where the right-hand side is defined to be zero on $\{f=0\}$ and belongs to $L^2$ by Lemma~\ref{lem:weighted-second}. The inequality also holds a.e. on $\{f=0\}$, since $\nabla f=0$ there. Thus the same lemma controls the singular gradient term near the zeros of $f$. For $p\ge1$ and $0<\varepsilon\le1$, the embedding $H^2\hookrightarrow W^{1,4}$ gives
\EQ{
f_\varepsilon^{p-1}|\nabla f|^2
\le(\normo{f}_\infty^2+1)^{(p-1)/2}|\nabla f|^2\in L^2.
}

The differentiated expressions converge pointwise almost everywhere on $\{f\ne0\}$. On $\{f=0\}$, they satisfy
\EQ{
D^2G_{\ell,\varepsilon}(f)
=e^{-i\lambda\ell\varepsilon^p}D^2f
\longrightarrow D^2f
\qquad\text{a.e.}
}
The preceding bounds therefore give convergence of the Hessians in $L^2$ by dominated convergence. Since $G_{\ell,\varepsilon}(f)\to G_\ell(f)$ in $L^2$, the limit is the distributional Hessian of $G_\ell(f)$, and hence
\EQ{
G_{\ell,\varepsilon}(f)\longrightarrow G_\ell(f)
\quad\text{in }H^2
\quad\text{as }\varepsilon\downarrow0.
}
Passing to the limit in the uniform Hessian estimate proves~\eqref{eq:phase-hessian}.
\end{proof}

The following standard estimates follow from the heat kernel representation and the Fourier estimates in~\cite[Sections~2.1.1--2.1.2]{BahouriCheminDanchin2011}, with heat time $1/t$.

\begin{lem}\label{lr:heat}
Let $W\in H^s(\R^d)$ with $d/2<s\le2$, and set $W_t:=e^{\Delta_y/t}W$ for $t\ge2$.
\begin{enumerate}
\item[(i)] For every $k\ge0$, the map $t\mapsto W_t$ belongs to $C^1([2,\infty);H^k)$, and
\EQn{\label{lr:heat-time-derivative}
\partial_tW_t=-t^{-2}\Delta_yW_t.
}
\item[(ii)] For every $t\ge2$,
\EQn{\label{lr:heat-bounds}
\normo{W_t}_2&\le\normo{W}_2, \\
 \normo{W_t}_\infty&\le\normo{W}_\infty,\\
 \normo{D^2W_t}_2+t^2\normo{\partial_tW_t}_2
 &\lsm_s t^{1-s/2}\normo{W}_{H^s}.
}
\item[(iii)] For every $t\ge2$ and $2\le r\le\infty$,
\EQn{\label{lr:heat-lr}
\normo{W_t-W}_r
 \lsm_{s,r} t^{-s/2+\frac d2(\frac12-\frac1r)}\normo{W}_{H^s}.
}
\end{enumerate}
\end{lem}

\subsection{A Fuchsian system}

The following system arises after freezing the spatial coefficients in the low-frequency argument of Subsection~\ref{lin:section-blocks}. For $g\in\R$, set
\EQn{\label{prelim:fuchsian-matrices}
A_g:=\begin{pmatrix}0&0\\g&0\end{pmatrix},\qquad
 J_0:=\begin{pmatrix}0&-1\\1&0\end{pmatrix},
}
and consider
\EQn{\label{prelim:fuchsian-system}
r\partial_r\vec{u}=(A_g+rJ_0)\vec{u},\qquad
 r>0,\quad\vec{u}(r)\in\C^2.
}

\begin{lem}[Fuchsian system with nilpotent residue]\label{prelim:fuchsian}
Let $G\ge0$ and $g\in\R$ satisfy $|g|\le G$. Every solution of~\eqref{prelim:fuchsian-system} satisfies
\EQn{\label{prelim:fuchsian-bound}
\normo{\vec{u}(r_1)}
 \lsm_G \brk{1+\abs{\log\frac{r_1}{r_0}}}
 \normo{\vec{u}(r_0)},
 \qquad r_0,r_1>0.
}
\end{lem}

\begin{proof}
Set $R_G:=2(G+1)$ and fix $0<r_0\le R_G$. Since $A_g^2=0$, variation of constants gives
\EQ{
\vec{u}(r)=\brk{I+A_g\log\frac r{r_0}}\vec{u}(r_0)
 +\int_{r_0}^r\brk{I+A_g\log\frac rs}J_0\vec{u}(s)\,\dd s,
 \qquad 0<r\le R_G.
}
For $s$ between $r$ and $r_0$, we have $|\log(r/s)|+|\log(s/r_0)|=|\log(r/r_0)|$, so
\EQ{
\frac{(1+G|\log(r/s)|)(1+G|\log(s/r_0)|)}
      {1+G|\log(r/r_0)|}
 &\le1+G\min\{|\log(r/s)|,|\log(s/r_0)|\}\\
 &\le1+G\log\frac{R_G}s.
}
Using $\normo{A_g}\le G$ and $\normo{J_0}=1$, we obtain
\EQ{
\frac{\normo{\vec{u}(r)}}{1+G|\log(r/r_0)|}
 \le\normo{\vec{u}(r_0)}
 +\abs{\int_{r_0}^r\brk{1+G\log\frac{R_G}s}
       \frac{\normo{\vec{u}(s)}}{1+G|\log(s/r_0)|}\,\dd s}.
}
Gronwall's inequality in either direction and
\EQ{
\int_0^{R_G}\brk{1+G\log\frac{R_G}s}\,\dd s=R_G(1+G)
}
give
\EQn{\label{lin:frozen-low-r}
\normo{\vec{u}(r_1)}
 \lsm_G
 \brk{1+\abs{\log\frac{r_1}{r_0}}}
 \normo{\vec{u}(r_0)},
 \qquad 0<r_0,r_1\le R_G.
}

For $r\ge R_G$, write $\vec{u}=(u_1,u_2)^T$ and set
\EQ{
E(r):=\brk{1+\frac gr}|u_1(r)|^2+|u_2(r)|^2.
}
Since $1/2\le1+g/r\le3/2$,
\EQ{
\frac12\normo{\vec{u}(r)}^2\le E(r)\le\frac32\normo{\vec{u}(r)}^2.
}
Equation~\eqref{prelim:fuchsian-system} gives $u_1'=-u_2$ and $u_2'=(1+g/r)u_1$. Hence
\EQn{\label{lin:frozen-high-r}
E'(r)=-\frac g{r^2}|u_1(r)|^2,
 \qquad |E'(r)|\le\frac{2G}{r^2}E(r).
}
For $r_0,r_1\ge R_G$, Gronwall's inequality applied to $E$ in either direction and $2G/R_G\le1$ give
\EQ{
\normo{\vec{u}(r_1)}^2 
 \le3\normo{\vec{u}(r_0)}^2\exp\brk{2
       \int_{\min\{r_0,r_1\}}^{\max\{r_0,r_1\}}
       \frac G{r^2}\,\dd r}
\lsm\normo{\vec{u}(r_0)}^2.
}
If $r_0$ and $r_1$ lie on opposite sides of $R_G$, applying the two estimates successively gives
\EQ{
\normo{\vec{u}(r_1)}
 \lsm_G\brk{1+\abs{\log\frac{r_1}{R_G}}
                  +\abs{\log\frac{R_G}{r_0}}}
                  \normo{\vec{u}(r_0)}
 =\brk{1+\abs{\log\frac{r_1}{r_0}}}
                  \normo{\vec{u}(r_0)}.
}
This proves~\eqref{prelim:fuchsian-bound}.
\end{proof}

\section{The linearized equation}\label{lin:section-frozen}

Let $W\in L^\infty(\R^d)\cap C^{0,\alpha}(\R^d)$ with $0<\alpha\le1$, and recall that $V(t)=G_{\log t}(W)$. Under $u=\mathcal P_tb$, with $\mathcal P_t$ defined in~\eqref{eq:ray-map}, equation~\eqref{eq:nls} becomes
\EQn{\label{eq:ray-equation}
i\partial_tb+\frac1{2t^2}\Delta_yb=\frac\lambda t|b|^pb.
}
Set
\EQ{
\gamma_0:=\min\{p,1\},\qquad \gamma:=\alpha\gamma_0.
}
The homogeneous linearized equation is for $t>0$,
\EQn{\label{lin:ray-equation}
i\partial_tz+\frac1{2t^2}\Delta_yz
 =\frac\lambda t\mbrkbb{
      \brk{1+\frac p2}|W|^p z +\frac p2|W|^{p-2}W^2e^{-2i\lambda|W|^p\log t}\wb z.
   }
}
For $W\in L^\infty$, both coefficients on the right-hand side of~\eqref{lin:ray-equation} belong to $C((0,\infty);L^\infty)$. The nonlinear estimate is
\EQ{
\normb{\frac\lambda tDF(V(t))[z]}_2
 \le\frac{1+p}{t}\normo{W}_\infty^p\normo{z}_2.
}
Then the existence of the solution to \eqref{lin:ray-equation} can be obtained by contraction mapping. Moreover, notice that the coefficient on the right of the nonlinear estimate is integrable on every compact subinterval of $(0,\infty)$. Then, since the contraction condition depends only on the time interval and $B$, successive applications of the contraction mapping argument give a unique solution on every compact subinterval of $(0,\infty)$. This yields the following result.

\begin{lem}[Solvability of the linearized Cauchy problem]\label{lin:cauchy-existence}
Let $W\in L^\infty(\R^d)$. For every $r>0$, equation~\eqref{lin:ray-equation} has a unique solution $z\in C((0,\infty);L^2)$ for each prescribed value $z(r)\in L^2$. 
\end{lem}

\subsection{Estimates on logarithmic time intervals}\label{lin:section-blocks}

Recall that $B=\normo{W}_\infty$. We now estimate solutions of~\eqref{lin:ray-equation} on $[e^\theta,e^{\theta+L}]$, an interval of logarithmic length $L$. Throughout this subsection, $\theta>1$ and $L\ge1$.

\begin{prop}\label{lin:block}
Let $0<\alpha\le1$ and $W\in L^\infty(\R^d)\cap C^{0,\alpha}(\R^d)$. There exists a constant $C_{p,B}\ge1$ such that, for every $L\ge1$, there exists $\Theta_L>1$ for which every solution $z$ of~\eqref{lin:ray-equation} satisfies
\EQn{\label{lin:block-estimate}
\normo{z(t)}_2\le C_{p,B}L\normo{z(r)}_2,
 \qquad e^\theta\le t\le r\le e^{\theta+L},\quad \theta\ge\Theta_L.
}
The threshold $\Theta_L$ may also depend on $L$, $\alpha$, and the H\"older seminorm of $W$.
\end{prop}

We reduce~\eqref{lin:block-estimate} to an estimate for the real and imaginary parts on $[0,L]$. Fix $L\ge1$ and $\theta>1$. To remove the factor $t^{-1}$ from the potential, set
\EQ{
t:=e^{\theta+\sigma},\qquad
 z_\theta(\sigma,y):=z(e^{\theta+\sigma},y),
 \qquad 0\le\sigma\le L.
}
Since $\partial_\sigma z_\theta=t\partial_tz$, equation~\eqref{lin:ray-equation} gives
\EQn{\label{lin:log-equation}
i\partial_\sigma z_\theta+\tfrac12e^{-\theta-\sigma}\Delta_yz_\theta
 =\lambda\brk{1+\frac p2}|W|^p z_\theta +\frac{\lambda p}{2}|W|^{p-2}W^2
 e^{-2i\lambda|W|^p(\theta+\sigma)}\wb{z_\theta}.
}
Set
\EQ{
X:=\begin{pmatrix}
 \operatorname{Re}z_\theta\\
 \operatorname{Im}z_\theta
 \end{pmatrix}\in L^2(\R^d)\times L^2(\R^d).
}
We denote the real Jacobian matrix of $F$ by
\EQ{
J_F(w):=|w|^pI+p|w|^{p-2}
 \begin{pmatrix}
 (\operatorname{Re}w)^2&\operatorname{Re}w\operatorname{Im}w\\
 \operatorname{Re}w\operatorname{Im}w&(\operatorname{Im}w)^2
 \end{pmatrix}.
}
For $w,h\in\C$, the matrix $J_F(w)$ represents $DF(w)$ through
\EQ{
\begin{pmatrix}
 \operatorname{Re}(DF(w)[h])\\
 \operatorname{Im}(DF(w)[h])
 \end{pmatrix}
 =J_F(w)\begin{pmatrix}
 \operatorname{Re}h\\
 \operatorname{Im}h
 \end{pmatrix}.
}
Multiplication by $i$ is represented by $J_0$ from~\eqref{prelim:fuchsian-matrices}. Taking real and imaginary parts of~\eqref{lin:log-equation} therefore gives
\EQ{
J_0\partial_\sigma X+\tfrac12e^{-\theta-\sigma}\Delta_yX
 =\lambda J_F\brkb{V(e^{\theta+\sigma})}X.
}
We separate the scalar part of this symmetric matrix. Define
\EQn{\label{lin:semiclassical-definitions}
D_\theta(\sigma):=-\tfrac12e^{-\theta-\sigma}\Delta_y,\qquad
 L_0(\sigma):=-D_\theta(\sigma)J_0,
}
and the coefficient matrices
\EQn{\label{lin:matrix-potentials}
\V&:=-\lambda\brk{1+\frac p2}|W|^pJ_0,\\
 \B(\sigma)&:=-\lambda J_0\mbrk{
 J_F\brkb{V(e^{\theta+\sigma})}
 -\brk{1+\frac p2}|W|^pI}.
}
Explicitly,
\EQ{
\B(\sigma)=\frac{\lambda p}{2}|W|^{p-2}
 \begin{pmatrix}
 \operatorname{Im}\!\brk{W^2e^{-2i\lambda|W|^p(\theta+\sigma)}}
 &-\operatorname{Re}\!\brk{W^2e^{-2i\lambda|W|^p(\theta+\sigma)}}\\
 -\operatorname{Re}\!\brk{W^2e^{-2i\lambda|W|^p(\theta+\sigma)}}
 &-\operatorname{Im}\!\brk{W^2e^{-2i\lambda|W|^p(\theta+\sigma)}}
 \end{pmatrix}.
}
Multiplying this real system by $-J_0$ gives
\EQn{\label{lin:doubled-equation}
\partial_\sigma X=(L_0+\V+\B)X.
}
Therefore, using this system, we reduce the proof of Proposition~\ref{lin:block} to the following: for every $L\ge1$, one can choose $\Theta_L>1$ such that every solution of~\eqref{lin:doubled-equation} satisfies
\EQ{
\normo{X(\sigma)}_2\lsm_{p,B} L\normo{X(\rho)}_2,
 \qquad 0\le\sigma\le\rho\le L,\quad\theta\ge\Theta_L.
}

Lemma~\ref{lin:coefficient-bounds} records the coefficient bounds used in the low- and high-frequency estimates of Lemmas~\ref{lin:low} and~\ref{lin:high}. We combine those estimates to prove Proposition~\ref{lin:block} after Lemma~\ref{lin:high}.

\begin{lem}[Coefficient bounds]\label{lin:coefficient-bounds}
Let $0<\alpha\le1$ and $W\in L^\infty(\R^d)\cap C^{0,\alpha}(\R^d)$. Fix $\theta>1$ and $L\ge1$, and let $\V,\B$ be defined by~\eqref{lin:matrix-potentials}. Then, for every $0\le\sigma\le L$,
\EQnn{
&\label{lin:matrix-structure}
 \V^*=-\V,\qquad \B^*=\B,
 \qquad J_0\B=-\B J_0,\\&
 \normo{\V}_\infty+\normo{\B(\sigma)}_\infty
 +\normo{\partial_\sigma\B(\sigma)}_\infty\lsm_{p,B} 1,
 \label{lin:matrix-bounded}\\&
 \normo{\B(\sigma)-\B(\rho)}_\infty
 \le pB^{2p}|\sigma-\rho|,\qquad 0\le\rho\le L,\label{lin:matrix-time-lipschitz}\\&
 [\V]_{C^{0,\gamma}}+[\B(\sigma)]_{C^{0,\gamma}}
 \lsm_{W} \theta+L.
 \label{lin:matrix-holder}
}
\end{lem}

\begin{proof}
We first prove the identities and bounds in~\eqref{lin:matrix-structure}--\eqref{lin:matrix-time-lipschitz}. The real matrix $J_F(w)$ is symmetric and satisfies
\EQ{
\operatorname{tr}J_F(w)=(p+2)|w|^p,\qquad
 \norm{J_F(w)-\brk{1+\frac p2}|w|^pI}
 =\frac p2|w|^p.
}
Every real symmetric $2\times2$ matrix of trace zero anticommutes with $J_0$. Since $|V(t)|=|W|$, this gives~\eqref{lin:matrix-structure} and
\EQ{
\normo{\V}_\infty=\brk{1+\frac p2}B^p,\qquad
 \normo{\B(\sigma)}_\infty=\frac p2 B^p.
}
Since $F(e^{i\varphi}w)=e^{i\varphi}F(w)$, differentiation in real coordinates gives
\EQ{
J_F(e^{i\varphi}w)e^{\varphi J_0}
 =e^{\varphi J_0}J_F(w).
}
Multiplying on the right by $e^{-\varphi J_0}$ yields the rotation identity
\EQ{
J_F(e^{i\varphi}w)=e^{\varphi J_0}J_F(w)e^{-\varphi J_0},
 \qquad \varphi\in\R.
}
Using this identity and $V(e^{\theta+\sigma})=e^{-i\lambda|W|^p(\theta+\sigma)}W$, we obtain
\EQ{
\partial_\sigma\B=-2\lambda|W|^pJ_0\B,\qquad
 \normo{\partial_\sigma\B}_\infty\le pB^{2p}.
}
This proves~\eqref{lin:matrix-bounded} and~\eqref{lin:matrix-time-lipschitz}.

We now prove~\eqref{lin:matrix-holder}. The maps $w\mapsto|w|^p$ and $w\mapsto J_F(w)$ are $\gamma_0$-H\"older continuous on bounded subsets of $\C$. By~\eqref{lr:phase-first-derivative}, for $\ell\ge0$ and $|w_1|,|w_2|\le B$,
\EQ{
|G_\ell(w_1)-G_\ell(w_2)|
 \le(1+p\ell B^p)|w_1-w_2|.
}
Taking $\ell=\theta+\sigma$ therefore gives
\EQ{
{}[\V]_{C^{0,\gamma}}+[\B(\sigma)]_{C^{0,\gamma}}
 &\lsm_{p,B} (\theta+L)^{\gamma_0}
       [W]_{C^{0,\alpha}}^{\gamma_0}\\
 &\lsm_{W} \theta+L,
}
since $\theta>1$, $\gamma=\alpha\gamma_0$, and $0<\gamma_0\le1$. This proves~\eqref{lin:matrix-holder}.
\end{proof}

For $K\ge1$, put $N:=Ke^{\theta/2}$ in the frequency cutoffs defined by~\eqref{lin:P-Q}. We write $\M:=\V+\B$.

Replacing $W$ by a constant $w_0\in\C$ in~\eqref{lin:doubled-equation}, we obtain the following initial value problem on $[0,L]$. Given $0\le\rho\le L$ and $X_0\in L^2(\R^d)\times L^2(\R^d)$, let $X_{w_0}$ solve
\EQn{\label{lin:frozen-equation}
\partial_\sigma X_{w_0}
 &=\tfrac12e^{-\theta-\sigma}J_0\Delta_yX_{w_0}
   -\lambda J_0J_F\!\brk{w_0e^{-i\lambda|w_0|^p(\theta+\sigma)}}X_{w_0},\\
 X_{w_0}(\rho)&=X_0.
}
\begin{samepage}
\begin{lem}[Frozen equation and low-frequency estimate]
\label{lin:frozen-ode}\label{lin:low}
Let $0<\alpha\le1$ and $W\in L^\infty(\R^d)\cap C^{0,\alpha}(\R^d)$. Fix $L,K\ge1$ and $\theta>1$, and set $N:=Ke^{\theta/2}$. Then the following estimates hold.
\begin{enumerate}
\item[(i)] For every $w_0\in\C$ with $|w_0|\le B$, the solution $X_{w_0}$ of~\eqref{lin:frozen-equation} satisfies
\EQn{\label{lin:frozen-pde-bound}
\normo{X_{w_0}(\sigma)}_2\lsm_{p,B} L\normo{X_{w_0}(\rho)}_2,
 \qquad 0\le\sigma,\rho\le L.
}

\item[(ii)] Let $X$ be a solution of~\eqref{lin:doubled-equation} and set $v:=P_{\le N}X$. Then for $0\le\sigma\le\rho\le L$,
\EQn{\label{lin:low-estimate}
\normo{v(\sigma)}_2 \le C_{p,B}L\normo{v(\rho)}_2  +C_{W,K}L(\theta+L)e^{-\gamma\theta/4}
       \int_\sigma^\rho\normo{X(r)}_2\,\dd r.
}
\end{enumerate}
\end{lem}
\end{samepage}

\begin{proof}
We first prove (i). Fix $0\le\rho\le L$ and let $X_0$ have real-valued Schwartz components. Let $X_{w_0}$ solve~\eqref{lin:frozen-equation} with $X_{w_0}(\rho)=X_0$.
Write $w_0=|w_0|e^{i\phi_0}$, taking $\phi_0:=0$ when $w_0=0$. To remove the phase of the frozen profile, define
\EQ{
Y(\sigma,y):=e^{(\lambda|w_0|^p(\theta+\sigma)-\phi_0)J_0}X_{w_0}(\sigma,y).
}
The rotation is orthogonal, so
\EQ{
\normo{Y(\sigma)}_2=\normo{X_{w_0}(\sigma)}_2,
 \qquad \normo{Y(\rho)}_2=\normo{X_0}_2.
}
Applying the rotation identity with $w=|w_0|$ and $\varphi=\phi_0-\lambda|w_0|^p(\theta+\sigma)$ gives
\EQ{
e^{(\lambda|w_0|^p(\theta+\sigma)-\phi_0)J_0}
 J_F\!\brk{w_0e^{-i\lambda|w_0|^p(\theta+\sigma)}}X_{w_0}
 =J_F(|w_0|)Y.
}
Since $J_F(|w_0|)=|w_0|^p\operatorname{diag}(1+p,1)$, differentiating $Y$ and using~\eqref{lin:frozen-equation} gives
\EQ{
\partial_\sigma Y
 &=\tfrac12e^{-\theta-\sigma}J_0\Delta_yY
   +\lambda|w_0|^pJ_0Y-\lambda J_0J_F(|w_0|)Y\\
 &=\tfrac12e^{-\theta-\sigma}J_0\Delta_yY-A_gY,
}
where $g:=\lambda p|w_0|^p$. Taking the Fourier transform in $y$ componentwise, we write $\wh Y(\sigma,\xi):=\mathcal F_y[Y(\sigma,\cdot)](\xi)$ and obtain
\EQn{\label{lin:frozen-ode-system}
\partial_\sigma\wh Y(\sigma,\xi)
 =-\brk{A_g+\tfrac12|\xi|^2e^{-\theta-\sigma}J_0}\wh Y(\sigma,\xi).
}
For a fixed spatial frequency $\xi\ne0$, set $r:=\tfrac12|\xi|^2e^{-\theta-\sigma}$ and regard $\wh Y$ as a function of $r$. Since $\partial_\sigma=-r\partial_r$, equation~\eqref{lin:frozen-ode-system} becomes
\EQ{
r\partial_r\wh Y=(A_g+rJ_0)\wh Y.
}
This is~\eqref{prelim:fuchsian-system}, so Lemma~\ref{prelim:fuchsian}, with $G=pB^p$, gives
\EQ{
|\wh Y(\sigma,\xi)|
  \lsm_{p,B} \brk{1+\abs{\log
       \frac{|\xi|^2e^{-\theta-\sigma}}{|\xi|^2e^{-\theta-\rho}}}}
       |\wh Y(\rho,\xi)| =(1+|\sigma-\rho|)|\wh Y(\rho,\xi)|.
}
For $\xi=0$, equation~\eqref{lin:frozen-ode-system} and $A_g^2=0$ give
\EQ{
\wh Y(\sigma,0)=\brkb{I-(\sigma-\rho)A_g}\wh Y(\rho,0).
}
Thus, for every $\xi\in\R^d$,
\EQn{\label{lin:frozen-ode-bound}
|\wh Y(\sigma,\xi)|
 \lsm_{p,B} (1+|\sigma-\rho|)|\wh Y(\rho,\xi)|.
}
By Plancherel's theorem and the norm identities above,
\EQ{
\normo{X_{w_0}(\sigma)}_2
 =\normo{Y(\sigma)}_2 \lsm_{p,B} (1+|\sigma-\rho|)\normo{Y(\rho)}_2 \lsm_{p,B} L\normo{X_0}_2.
}
Approximation of $X_0$ in $L^2$ proves~\eqref{lin:frozen-pde-bound} for arbitrary data.

We now prove (ii) by spatial localization and (i). Since $P_{\le N}$ commutes with $L_0$, the full equation for $v$ is
\EQ{
\partial_\sigma v=(L_0+\M )v+[P_{\le N},\M ]X.
}
Choose smooth real cutoff functions $\{\chi_\nu\}_{\nu\in\mathbb Z^d}$ at scale $e^{-\theta/4}$, with bounded overlap and satisfying
\EQn{\label{lin:square-partition}
\sum_\nu\chi_\nu^2=1,\qquad
 \sum_\nu|\nabla\chi_\nu|^2\lsm e^{\theta/2},\qquad
 \sum_\nu|\Delta\chi_\nu|^2\lsm e^\theta.
}
Each $\chi_\nu$ is supported in a ball of radius $Ce^{-\theta/4}$ centered at $y_\nu$. These functions are obtained from lattice translates of a fixed bump covering $\mathbb R^d$, after rescaling and division by the square root of the sum of their squares.

Set $v_\nu:=\chi_\nu v$ and define $\M _\nu$ by replacing $W(y)$ with $W(y_\nu)$ in \eqref{lin:matrix-potentials}. Then
\EQn{\label{lin:localized-equation}
\partial_\sigma v_\nu=(L_0+\M _\nu)v_\nu+f_\nu,
}
with
\EQn{\label{lin:localized-forcing}
f_\nu:=\chi_\nu(\M -\M _\nu)v
       -[L_0,\chi_\nu]v+\chi_\nu[P_{\le N},\M ]X.
}
For the coefficient error, \eqref{lin:matrix-holder}, the support condition and \eqref{lin:square-partition} give
\EQ{
\brk{\sum_\nu
 \normo{\chi_\nu(\M -\M _\nu)v}_2^2}^{1/2}
 \lsm_{W} (\theta+L)e^{-\gamma\theta/4}\normo{v}_2.
}
The $\ell^2_\nu L^2$ norm of the last term in \eqref{lin:localized-forcing} equals $\normo{[P_{\le N},\M ]X}_2$. Estimates~\eqref{lin:P-commutator} and~\eqref{lin:matrix-holder}, with $N=Ke^{\theta/2}$, give
\EQ{
\normo{[P_{\le N},\M]X}_2
 \lsm_{W} (\theta+L)K^{-\gamma}e^{-\gamma\theta/2}\normo{X}_2 \lsm_{W} (\theta+L)e^{-\gamma\theta/4}\normo{X}_2,
}
since $K\ge1$ and $\theta>1$.

For the Laplacian commutator, the product rule gives
\EQ{
[L_0,\chi_\nu]v
 =\tfrac12e^{-\theta-\sigma}J_0
 \brkb{2\nabla\chi_\nu\cdot\nabla v+(\Delta\chi_\nu)v}.
}
Consequently,
\EQnn{\label{lin:laplacian-commutator}
\brk{\sum_\nu\normo{[L_0,\chi_\nu]v}_2^2}^{1/2}
 \lsm e^{-\theta}
       \brk{e^{\theta/4}\normo{\nabla v}_2+
             e^{\theta/2}\normo{v}_2} \lsm Ke^{-\gamma\theta/4}\normo{X}_2.
}
In the second line, the frequency support of $v=P_{\le N}X$ gives $\normo{\nabla v}_2\le2Ke^{\theta/2}\normo{X}_2$. The last line uses $K\ge1$ and $0<\gamma\le1$. Summing the three bounds yields
\EQn{\label{lin:forcing-square-sum}
\brk{\sum_\nu\normo{f_\nu(r)}_2^2}^{1/2}
 \lsm_{W,K} (\theta+L)e^{-\gamma\theta/4}\normo{X(r)}_2.
}
Apply Duhamel's formula for the frozen equation to \eqref{lin:localized-equation}, with data at $\rho$, and estimate each homogeneous solution by~\eqref{lin:frozen-pde-bound}. Taking the norm in $\ell^2_\nu L^2$, Minkowski's inequality and \eqref{lin:forcing-square-sum} give
\EQ{
\brkb{\sum_\nu\normo{v_\nu(\sigma)}_2^2}^{\frac12}
 \le C_{p,B}L\brkb{\sum_\nu\normo{v_\nu(\rho)}_2^2}^{\frac12} +C_{W,K}L(\theta+L)e^{-\frac{\ga\theta}{4}}
       \int_\sigma^\rho\normo{X(r)}_2\,\dd r.
}
This proves \eqref{lin:low-estimate}.
\end{proof}

\begin{samepage}
\begin{lem}[High-frequency estimate]\label{lin:high}
Let $0<\alpha\le1$ and $W\in L^\infty(\R^d)\cap C^{0,\alpha}(\R^d)$. Fix $L\ge1$. Then one can choose $K\ge1$ such that the following holds for every $\theta>1$. If $X$ solves~\eqref{lin:doubled-equation}, set $N:=Ke^{\theta/2}$ and $w:=P_{>N}X$. Then
\EQn{\label{lin:high-estimate}
\normo{w(\sigma)}_2 \le2\normo{w(\rho)}_2 +C_W(\theta+L)e^{-\gamma\theta/4}
       \int_\sigma^\rho\normo{X(r)}_2\,\dd r
}
for $0\le\sigma\le\rho\le L$.
\end{lem}
\end{samepage}

\begin{proof}
We use a modified $L^2$ energy to cancel the leading growth caused by the symmetric coefficient $\B$. With $\ph_\theta$ as in~\eqref{lin:mollifier}, set
\EQ{
\V_\theta:=\ph_\theta*\V,\qquad
 \B_\theta(\sigma):=\ph_\theta*\B(\sigma).
}
Since $\ph_\theta$ is real, these coefficients satisfy the three identities in~\eqref{lin:matrix-structure}. The convolution bounds~\eqref{lin:mollification-estimates} and the coefficient bound~\eqref{lin:matrix-bounded} give
\EQnn{
&\label{lin:mollified-time-bounds}
 \normo{\V _\theta}_\infty+\normo{\B _\theta}_\infty
 +\normo{\partial_\sigma\B _\theta}_\infty\lsm_{p,B} 1,\\&
 \normo{\nabla\B _\theta}_\infty\lsm_{p,B} e^{\theta/4},\qquad
 \normo{\Delta\B _\theta}_\infty\lsm_{p,B} e^{\theta/2}.
 \label{lin:mollified-space-bounds}
}
Writing $\M _\theta:=\V _\theta+\B _\theta$ and $\M w=\M_\theta w+(\M-\M_\theta)w$, we obtain
\EQn{\label{lin:high-equation}
\partial_\sigma w=(L_0+\M _\theta)w+f_h,
}
with
\EQn{\label{lin:high-forcing}
f_h:=(\M-\M_\theta)w+[P_{>N},\M]X.
}
By~\eqref{lin:P-commutator}, \eqref{lin:mollification-estimates}, and~\eqref{lin:matrix-holder}, using $\normo{w}_2\le\normo{X}_2$, we have
\EQ{
\normo{f_h}_2
\le\normo{(\M-\M_\theta)w}_2
+\normo{[P_{>N},\M]X}_2 \lsm_{W} (\theta+L)e^{-\gamma\theta/4}\normo{X}_2.
}

Set
\EQn{\label{lin:R-definition}
\mathcal R(\sigma):=D_\theta(\sigma)^{-1}P_{>N/2}.
}
Define the symbol to be zero near the origin. The multiplier bounds are
\EQnn{
&\label{lin:R-bounds}
 \normo{\mathcal R(\sigma)u}_2\lsm \frac{e^\sigma}{K^2}\normo{u}_2,\qquad
 \normo{\nabla\mathcal R(\sigma)u}_2
       \lsm \frac{e^{\sigma+\theta/2}}{K}\normo{u}_2,\qquad
 \mathcal R'=\mathcal R,\\&
 \mathcal R D_\theta w=D_\theta\mathcal Rw=w.
 \label{lin:R-inverse-on-w}
}
Here $u$ is any vector in $L^2(\R^d)\times L^2(\R^d)$. 

Set
\EQn{\label{lin:K-definition}
\K:=\tfrac12J_0\B_\theta\mathcal R,
 \qquad E(\sigma):=\normo{w(\sigma)}_2^2
       +2\jb{\K(\sigma)w(\sigma),w(\sigma)},
}
where $\jb{\cdot,\cdot}$ denotes the real $L^2$ inner product. Since $J_0\B_\theta J_0=\B_\theta$ and $D_\theta,\mathcal R,J_0$ commute, we have
\EQnn{\label{lin:exact-commutator}
[\K,L_0]w
 = -\tfrac12\B_\theta\mathcal RD_\theta w
    -\tfrac12D_\theta\B_\theta\mathcal Rw = -\B_\theta w
    -\tfrac12[D_\theta,\B_\theta]\mathcal Rw.
}
We compute $E'$ using~\eqref{lin:high-equation}. The product rule gives
\EQ{
\frac12E'
 =\jb{\partial_\sigma w,w}+\jb{\K'w,w}
  +\jb{\K\partial_\sigma w,w}
  +\jb{\K w,\partial_\sigma w}.
}
Since $(L_0+\V_\theta)^*=-(L_0+\V_\theta)$, the first term is
\EQ{
\jb{\partial_\sigma w,w}
 =\jb{\B_\theta w,w}+\jb{f_h,w}.
}
For the last two terms, the contributions of $L_0+\V_\theta$ satisfy
\EQ{
 \jb{\K(L_0+\V_\theta)w,w}
   +\jb{\K w,(L_0+\V_\theta)w} 
 & =\jb{\K(L_0+\V_\theta)w,w}
          -\jb{(L_0+\V_\theta)\K w,w}\\
 & =\jb{\brkb{[\K,L_0]+[\K,\V_\theta]}w,w}.
}
Similarly, $\B_\theta^*=\B_\theta$ gives
\EQ{
 \jb{\K\B_\theta w,w}+\jb{\K w,\B_\theta w} =\jb{\K\B_\theta w,w}+\jb{\B_\theta\K w,w} =\jb{(\K\B_\theta+\B_\theta\K)w,w}.
}
The three forcing terms combine as
\EQ{
 \jb{f_h,w}+\jb{\K f_h,w}+\jb{\K w,f_h} 
 =\jb{f_h,(I+\K+\K^*)w}.
}
Combining these identities and then applying~\eqref{lin:exact-commutator}, we obtain
\EQnn{\label{lin:energy-identity}
\frac12E'
 &=\jb{\brkb{\B_\theta+[\K,L_0]}w,w}
   +\jb{\K'w,w} +\jb{\brkb{[\K,\V_\theta]
                  +\K\B_\theta+\B_\theta\K}w,w}\notag\\
 &\quad+\jb{f_h,(I+\K+\K^*)w}\notag\\
 &=-\tfrac12\jb{[D_\theta,\B_\theta]\mathcal Rw,w}
   +\jb{\K'w,w} +\jb{[\K,\V_\theta]w,w}
   +\jb{\K\B_\theta w,w}\notag\\
 &\quad+\jb{\B_\theta\K w,w}+\jb{f_h,(I+\K+\K^*)w}.
}

By~\eqref{lin:mollified-time-bounds}, \eqref{lin:R-bounds}, and $\mathcal R'=\mathcal R$, every $u\in L^2(\R^d)\times L^2(\R^d)$ satisfies
\EQn{\label{lin:K-bounds}
\normo{\K u}_2+\normo{\K'u}_2
 \lsm_{p,B}\frac{e^L}{K^2}\normo{u}_2.
}
Indeed, $\K'=\tfrac12J_0(\partial_\sigma\B_\theta)\mathcal R+\K$. Multiplication by $\V_\theta$ and $\B_\theta$ is bounded by~\eqref{lin:mollified-time-bounds}, so the three remaining products satisfy
\EQ{
\normo{[\K,\V_\theta]u}_2
 +\normo{\K\B_\theta u}_2+\normo{\B_\theta\K u}_2
 \lsm_{p,B}\frac{e^L}{K^2}\normo{u}_2.
}
For the commutator with $D_\theta$, the product rule gives the exact identity
\EQ{
[D_\theta,\B_\theta]\mathcal Ru
 =-\tfrac12e^{-\theta-\sigma}
   \brkb{2\nabla\B_\theta\cdot\nabla\mathcal Ru
                  +(\Delta\B_\theta)\mathcal Ru}.
}
Using~\eqref{lin:mollified-space-bounds} and~\eqref{lin:R-bounds}, we obtain
\EQnn{\label{lin:high-commutator-bound}
\normo{[D_\theta,\B_\theta]\mathcal Ru}_2
 \lsm_{p,B}e^{-\theta-\sigma}\brk{
       \frac{e^{\sigma+3\theta/4}}K
       +\frac{e^{\sigma+\theta/2}}{K^2}}\normo{u}_2 \lsm_{p,B}\frac1K\normo{u}_2.
}
In particular, each of the five quadratic terms in the last three lines of~\eqref{lin:energy-identity} is bounded by a constant times $e^L K^{-1}\normo{w}_2^2$. Thus
\EQn{\label{lin:energy-differential-bound}
\frac12|E'|
 \le C_{p,B}\frac{e^L}{K}\normo{w}_2^2
       +\normo{f_h}_2\brkb{\normo{w}_2+\normo{\K w}_2+\normo{\K^*w}_2}.
}
The same bound as in~\eqref{lin:K-bounds} holds for $\K^*u$, by the explicit formula $\K^*=\tfrac12\mathcal R J_0\B_\theta$.

For $w\in C([0,L];L^2)$, the integrated form of~\eqref{lin:energy-identity} follows by applying the calculation to $P_{\le m}w$, including the corresponding cutoff commutator in the forcing, and letting $m\to\infty$.

For fixed $L$, choose $K$ so that $Le^L/K$ is sufficiently small. The preceding bounds give, for every $u\in L^2(\R^d)\times L^2(\R^d)$,
\EQ{
\normo{\K u}_2\le\tfrac18\normo{u}_2,\qquad
 \normo{\K^*u}_2\le\tfrac18\normo{u}_2,\qquad
 C_{p,B}\frac{e^L}{K}\le\frac1{8L}.
}
The modified energy then satisfies
\EQn{\label{lin:energy-equivalence}
\frac34\normo{w(\sigma)}_2^2\le E(\sigma)
 \le\frac54\normo{w(\sigma)}_2^2.
}
By~\eqref{lin:energy-differential-bound},
\EQ{
\frac12|E'|
 \le\frac1{8L}\normo{w}_2^2
       +\frac54\normo{f_h}_2\normo{w}_2 \le\frac1{6L}E+\frac5{2\sqrt3}\normo{f_h}_2E^{1/2}.
}
Gronwall's inequality applied backward in time to $E^{1/2}$ gives
\EQ{
E(\sigma)^{1/2}
 \le e^{1/6}E(\rho)^{1/2}
       +\frac5{2\sqrt3}e^{1/6}
          \int_\sigma^\rho\normo{f_h(r)}_2\,\dd r.
}
Using~\eqref{lin:energy-equivalence} and the bound for $f_h$, we obtain
\EQ{
\normo{w(\sigma)}_2
 \le\sqrt{\frac53}e^{1/6}\normo{w(\rho)}_2
   +C_W(\theta+L)e^{-\gamma\theta/4}
       \int_\sigma^\rho\normo{X(r)}_2\,\dd r.
}
Since $\sqrt{5/3}\,e^{1/6}<2$, this proves~\eqref{lin:high-estimate}.
\end{proof}

\begin{proof}[Proof of Proposition~\ref{lin:block}]
Fix $L$ and choose $K$ as in Lemma~\ref{lin:high}. For a solution $X$ of~\eqref{lin:doubled-equation}, write $v:=P_{\le N}X$ and $w:=P_{>N}X$. Since $X=v+w$ and $\normo{v(\rho)}_2,\normo{w(\rho)}_2\le\normo{X(\rho)}_2$, Lemma~\ref{lin:low}(ii) and Lemma~\ref{lin:high} give
\EQn{\label{lin:combined-block}
\normo{X(\sigma)}_2
 \le C_{p,B}L\normo{X(\rho)}_2 +C_{W,K}L(\theta+L)e^{-\gamma\theta/4}
   \int_\sigma^\rho\normo{X(r)}_2\,\dd r.
}
For these fixed $L,K$, choose $\Theta_L>1$ large enough that
\EQ{
C_{W,K}L^2(\theta+L)e^{-\gamma\theta/4}\le\log2,
 \qquad\theta\ge\Theta_L.
}
This is possible because $\gamma>0$. Gronwall's inequality backward in time gives
\EQ{
\normo{X(\sigma)}_2
 \lsm_{p,B} L\normo{X(\rho)}_2.
}
The equality $\normo{X(\sigma)}_2=\normo{z(e^{\theta+\sigma})}_2$ proves~\eqref{lin:block-estimate} after increasing $C_{p,B}$.
\end{proof}

\subsection{An arbitrarily small power loss}\label{lin:subsec-iteration}

\begin{prop}\label{lin:main}
Let $0<\alpha\le1$ and suppose that $W\in L^\infty(\R^d)\cap C^{0,\alpha}(\R^d)$. For every $\eta>0$, there exists $T_{W,\eta}\ge1$ such that every solution $z$ of~\eqref{lin:ray-equation} satisfies
\EQn{\label{lin:main-estimate}
\normo{z(t)}_2
 \lsm_{W,\eta}\brk{\frac rt}^\eta\normo{z(r)}_2,
 \qquad r\ge t\ge T_{W,\eta}.
}
We impose no smallness condition on $W$.
\end{prop}

\begin{proof}
Given $\eta>0$, fix the constant $C_{p,B}\ge1$ in Proposition~\ref{lin:block} and choose $L\ge1$ so that
\EQn{\label{lin:L-choice}
\frac{\log[C_{p,B}L]}{L}<\eta.
}
Use the threshold $\Theta_L$ from Proposition~\ref{lin:block} for this $L$.

For $\tau_2\ge\tau_1\ge\Theta_L$, divide $[\tau_1,\tau_2]$ into intervals of length $L$ and at most one shorter interval. There are at most $1+(\tau_2-\tau_1)/L$ intervals. Estimate~\eqref{lin:block-estimate} applies to the shorter interval as well, since it holds for every pair of points in a full interval. Applying that estimate successively to the same solution gives
\EQ{
\normo{z(e^{\tau_1})}_2
 \le\mbrkb{C_{p,B}L}^{1+(\tau_2-\tau_1)/L}\normo{z(e^{\tau_2})}_2 \lsm_{W,\eta}e^{\eta(\tau_2-\tau_1)}\normo{z(e^{\tau_2})}_2.
}
Taking $T_{W,\eta}:=e^{\Theta_L}$ proves~\eqref{lin:main-estimate} after returning to $t=e^{\tau_1}$ and $r=e^{\tau_2}$.
\end{proof}

We next return to the original variables $(t,x)$. Recall that $A_W(t)=\mathcal P_tV(t)$. If $z$ solves
\EQ{
i\partial_tz+\frac1{2t^2}\Delta_yz
 =\frac\lambda tDF(V(t))[z]+\mathcal P_t^{-1}f,
}
then $h:=\mathcal P_tz$ solves
\EQn{\label{lin:physical-equation}
\brk{i\partial_t+\frac12\Delta_x}h
 =\lambda\brk{1+\frac p2}|A_W|^ph
   +\frac{\lambda p}{2}|A_W|^{p-2}A_W^2\wb h+f.
}

\begin{cor}[Estimate for the forced linearized equation]\label{lin:forced-estimate}
Let $0<\alpha\le1$ and $W\in L^\infty(\R^d)\cap C^{0,\alpha}(\R^d)$. For every $\eta>0$, there exists $T_{W,\eta}\ge2$ such that, for every $r\ge T_{W,\eta}$ and $f\in L^1([T_{W,\eta},r];L^2)$, each solution $h\in C([T_{W,\eta},r];L^2)$ of~\eqref{lin:physical-equation} satisfies
\EQn{\label{lin:physical-estimate}
\normo{h(t)}_2
 \lsm_{W,\eta}
 \brk{\frac rt}^\eta\normo{h(r)}_2
 +\int_t^r\brk{\frac\tau t}^\eta\normo{f(\tau)}_2\,\dd \tau,
 \qquad T_{W,\eta}\le t\le r.
}
\end{cor}

\begin{proof}
Proposition~\ref{lin:main} and the unitarity of $\mathcal P_t$ give~\eqref{lin:physical-estimate} for $f=0$. The general case follows from Duhamel's formula.
\end{proof}

\section{Heat regularization and the asymptotic profile}\label{sec:profile}

Let $W\in H^s(\R^d)$ with $d/2<s\le2$, and fix $\alpha:=\frac12\min\fbrk{1,s-\frac d2}>0$. 
Sobolev embedding gives
\EQ{
\normo{W}_\infty+[W]_{C^{0,\alpha}}\lsm_s\normo{W}_{H^s}.
}
Thus Corollary~\ref{lin:forced-estimate} applies to $W$ with this choice of $\alpha$.

We then find $u=\wt A+h$ with
\EQn{\label{lr:profiles}
W_t:=e^{\Delta_y/t}W,\qquad
 \wt V(t):=G_{\log t}(W_t),\qquad
 \wt A(t):=\mathcal P_t\wt V(t).
}
By Lemma~\ref{lr:heat}, $W_t\in H^2$ and
\EQ{
t^{-2}\normo{D^2W_t}_2+\normo{\partial_tW_t}_2
 \lsm_s t^{-1-s/2}\normo{W}_{H^s}.
}

\subsection{Error estimate for the regularized profile}

\begin{prop}\label{prop:error}
Let $W\in H^s(\R^d)$ with $d/2<s\le2$. Let $\wt A$ be the regularized profile defined in~\eqref{lr:profiles}. Then $\wt A\in C([2,\infty);L^2)$ and satisfies
\EQn{\label{eq:profile-equation}
\brk{i\partial_t+\frac12\Delta}\wt A
 =\lambda|\wt A|^p\wt A+\wt e,
}
where the error is given by
\EQn{\label{eq:error-exact}
\wt e(t):=\mathcal P_t\brk{
    iDG_{\log t}(W_t)[\partial_tW_t]
    +\frac1{2t^2}\Delta\wt V(t)}.
}
For every $t\ge2$, the error and the profile satisfy
\EQn{\label{eq:error-size}
\normo{\wt e(t)}_2\lsm_{W,s} t^{-1-s/2}(\log t)^2,
 \qquad
 \normo{\wt A(t)}_2=\normo{W_t}_2\le\normo{W}_2,
}
and
\EQn{\label{lr:auxiliary-pointwise}
\normo{\wt A(t)}_\infty^p
 =t^{-1}\normo{W_t}_\infty^p\le t^{-1}B^p.
}
Moreover,
\EQ{
\wt A\in L^q_{\mathrm{loc}}([2,\infty);L^r),\qquad
 \frac2q+\frac dr=\frac d2,\quad
 2<q\le\infty,\quad 2\le r\le\infty.
}
\end{prop}

\begin{proof}
Since $dp/2=1$, we have
\EQn{\label{eq:ray-identities}
\brk{i\partial_t+\frac12\Delta_x}\wt A-\lambda|\wt A|^p\wt A =\mathcal P_t\brk{
 i\partial_t\wt V-\frac\lambda t|\wt V|^p\wt V
 +\frac1{2t^2}\Delta_y\wt V}.
}
By Lemma~\ref{lr:heat}, $W_t\in C^1_{\rm loc}(L^2\cap L^\infty)$. The chain rule gives
\EQn{\label{lr:auxiliary-time-equation}
i\partial_t\wt V
=\frac\lambda t|\wt V|^p\wt V
   +iDG_{\log t}(W_t)[\partial_tW_t].
}
Combining~\eqref{eq:ray-identities} and~\eqref{lr:auxiliary-time-equation} gives~\eqref{eq:profile-equation} with the error~\eqref{eq:error-exact}.

By~\eqref{lr:phase-first-derivative} and~\eqref{lr:heat-bounds},
\EQn{\label{lr:time-error-bound}
\normo{DG_{\log t}(W_t)[\partial_tW_t]}_2
\le(1+pB^p\log t)\normo{\partial_tW_t}_2 \lsm_{W,s} t^{-1-s/2}\log t.
}
Estimates~\eqref{eq:phase-hessian} and~\eqref{lr:heat-bounds} give
\EQn{\label{lr:auxiliary-hessian}
\normo{\Delta\wt V(t)}_2
 \lsm_p (1+B^p\log t)^2\normo{D^2W_t}_2 \lsm_{W,s} t^{1-s/2}(\log t)^2.
}
Consequently, since $\normo{\mathcal P_tf}_2=\normo{f}_2$,
\EQ{
\normo{\wt e(t)}_2
 \le\normo{DG_{\log t}(W_t)[\partial_tW_t]}_2
+\frac1{2t^2}\normo{\Delta\wt V(t)}_2 \lsm_{W,s} t^{-1-s/2}(\log t)^2.
}
\eqref{lr:auxiliary-time-equation} and the preceding bounds give
\EQ{
\wt V\in C^1_{\rm loc}([2,\infty);L^2)
\cap L^\infty_{\rm loc}([2,\infty);H^2),
}
which justifies~\eqref{eq:ray-identities}.

Since $|G_\ell(z)|=|z|$,
\EQ{
\normo{\wt A(t)}_2&=\normo{\wt V(t)}_2=\normo{W_t}_2\le\normo{W}_2,\\
\normo{\wt A(t)}_\infty^p
&=t^{-dp/2}\normo{\wt V(t)}_\infty^p
=t^{-1}\normo{W_t}_\infty^p\le t^{-1}B^p.
}
The explicit formula for $\mathcal P_t$ and density in $L^2$ show that $t\mapsto\mathcal P_tf$ is continuous for each $f\in L^2$. Thus, as $t\to r\ge2$,
\EQ{
\normo{\wt A(t)-\wt A(r)}_2
\le\normo{\wt V(t)-\wt V(r)}_2
+\normo{(\mathcal P_t-\mathcal P_r)\wt V(r)}_2\longrightarrow0.
}

Finally,
\EQ{
\normo{|\wt A(t)|^p\wt A(t)}_2
\le\normo{\wt A(t)}_\infty^p\normo{\wt A(t)}_2
\le t^{-1}B^p\normo{W}_2.
}
Hence $\lambda|\wt A|^p\wt A+\wt e\in L^1_{\rm loc}([2,\infty);L^2)$. Since $\wt A\in C([2,\infty);L^2)$, the distributional identity~\eqref{eq:profile-equation} implies the Duhamel formula for $t,t_0\ge2$,
\EQn{\label{lr:auxiliary-mild}
\wt A(t)
 =e^{i(t-t_0)\Delta/2}\wt A(t_0)-i\int_{t_0}^t e^{i(t-\sigma)\Delta/2}
       \brkb{\lambda|\wt A(\sigma)|^p\wt A(\sigma)
                        +\wt e(\sigma)}\,\dd \sigma.
}
The stated spacetime estimate follows from the Strichartz estimates.
\end{proof}

\subsection{Comparison with the prescribed profile}

\begin{prop}\label{lr:profile-comparison}
Let $W\in H^s(\R^d)$ with $d/2<s\le2$. Let $A_W$ and $\wt A$ be defined by~\eqref{eq:profile} and~\eqref{lr:profiles}. Then $A_W$ belongs to $C([2,\infty);L^2)\cap L^{\frac{2(d+2)}d}_{\mathrm{loc}}([2,\infty);L^{\frac{2(d+2)}d})$, and for every $R\ge2$,
\EQn{\label{lr:profile-comparison-estimate}
\normo{\wt A-A_W}_{L_t^\infty L_x^2\cap L_{t,x}^{\frac{2(d+2)}d}([R,2R])} \lsm_{W,s} R^{-s/2}\log R.
}
\end{prop}

\begin{proof}
To prove~\eqref{lr:profile-comparison-estimate}, we estimate the identity
\EQ{
 \wt A(t)-A_W(t)
 =\mathcal P_t\brkb{G_{\log t}(W_t)-G_{\log t}(W)}
}
in $L^2$ and $L^{\frac{2(d+2)}d}$. The bound on $DG_\ell$ in~\eqref{lr:phase-first-derivative} and the real mean value formula on the disk $|z|\le B$ give
\EQn{\label{lr:phase-lipschitz}
|G_{\log t}(z)-G_{\log t}(w)|
 \lsm_{p,B}\log t\,|z-w|.
}
Using $\normo{W_t}_\infty\le B$, Lemma~\ref{lr:heat}, and $\normo{\mathcal P_tf}_2=\normo{f}_2$, we obtain
\EQn{\label{lr:profile-comparison-l2}
\normo{\wt A(t)-A_W(t)}_2
 \lsm_{W,s} t^{-s/2}\log t.
}
For the spacetime part, admissibility gives
\EQ{
 \normo{\mathcal P_tf}_{\frac{2(d+2)}d}=t^{-d/(d+2)}\normo{f}_{\frac{2(d+2)}d}.
}
Thus \eqref{lr:heat-lr} and \eqref{lr:phase-lipschitz} give the spatial $L^{\frac{2(d+2)}d}$ bound
\EQn{\label{lr:profile-comparison-lq}
\normo{\wt A(t)-A_W(t)}_{\frac{2(d+2)}d}
 \lsm_{W,s} t^{-d/(d+2)}\log t\,\normo{W_t-W}_{\frac{2(d+2)}d}
 \lsm_{W,s} t^{-s/2-d/(2(d+2))}\log t.
}
Integrating in $t$ over $[R,2R]$, we obtain
\EQ{
\normo{\wt A-A_W}_{L^{\frac{2(d+2)}d}_{t,x}([R,2R])}^{\frac{2(d+2)}d}
 \lsm_{W,s}R^{-\frac{s(d+2)}d}(\log R)^{\frac{2(d+2)}d}.
}
Combining this with~\eqref{lr:profile-comparison-l2} proves~\eqref{lr:profile-comparison-estimate}. The formula~\eqref{eq:profile} gives $A_W\in C([2,\infty);L^2)$. The spacetime estimate follows by Proposition~\ref{prop:error} and~\eqref{lr:profile-comparison-estimate}
\end{proof}

\subsection{The equation linearized around the regularized profile}\label{lin:subsec-perturbation}

We show that the growth estimate of Corollary~\ref{lin:forced-estimate} persists for the equation linearized around $\wt A$.

\begin{prop}[The regularized linearized equation]\label{lr:perturbed-flow}
Let $W\in H^s(\R^d)$ with $d/2<s\le2$, and let $\wt A$ be defined by~\eqref{lr:profiles}. For every $r\ge2$ and $f\in L^1([2,r];L^2)$, the equation
\EQn{\label{lr:regularized-linear-equation}
\brk{i\partial_t+\frac12\Delta}h
 =\lambda DF(\wt A(t))[h]+f
}
has a unique solution $h\in C([2,r];L^2)$ for each prescribed value $h(r)\in L^2(\R^d)$. For every $\eta>0$, there exists $T_{W,s,\eta}\ge2$ such that
\EQn{\label{lr:perturbed-flow-bound}
\normo{h(t)}_2
 \lsm_{W,s,\eta}
       \brk{\frac rt}^\eta\normo{h(r)}_2
       +\int_t^r\brk{\frac\tau t}^\eta\normo{f(\tau)}_2\,\dd \tau,
 \qquad r\ge t\ge T_{W,s,\eta}.
}
\end{prop}

\begin{proof}
By~\eqref{lr:auxiliary-pointwise},
\EQ{
|DF(\wt A(t,x))[z]|
 \le(1+p)|\wt A(t,x)|^p|z|
 \le(1+p)B^pt^{-1}|z|.
}
We solve the free Duhamel formula by the contraction mapping:
\EQn{\label{lr:perturbed-duhamel}
h(t)
 =e^{i(t-r)\Delta/2}h(r) +i\int_t^r e^{i(t-\tau)\Delta/2}
       \brkb{\lambda DF(\wt A(\tau))[h(\tau)]+f(\tau)}\,\dd \tau.
}
For $I:=[a,r]$ with $2\le a<r$, unitarity of the free evolution gives
\EQ{
 \norm{\int_t^r e^{i(t-\tau)\Delta/2}
       DF(\wt A(\tau))[h(\tau)]\,\dd \tau}_{L^\infty(I;L^2)} \le(1+p)B^p\log\frac ra\,\normo{h}_{L^\infty(I;L^2)}.
}
Since the potential term is real-linear in $h$, the right-hand side of~\eqref{lr:perturbed-duhamel}, with $h(r)$ fixed, is a contraction on $C(I;L^2)$ whenever $(1+p)B^p\log(r/a)\le1/2$. Applying this argument on finitely many consecutive intervals gives a unique solution on $[2,r]$. The Strichartz estimates, applied to~\eqref{lr:perturbed-duhamel}, give
\EQ{
h\in L^{q_1}([2,r];L^{r_1}),\qquad
 \frac2{q_1}+\frac d{r_1}=\frac d2,\quad
 2<q_1\le\infty,\quad 2\le r_1\le\infty.
}

We now prove~\eqref{lr:perturbed-flow-bound}. Set
\EQ{
\gamma_0:=\min\{p,1\},\qquad
 \nu:=\gamma_0\brk{\frac s2-\frac d4}>0.
}
By Lemma~\ref{lr:heat},
\EQ{
\normo{W_t}_\infty\le B,\qquad
\normo{W_t-W}_\infty\lsm_s t^{d/4-s/2}\normo{W}_{H^s}.
}
Applying~\eqref{lr:phase-lipschitz} pointwise yields
\EQn{\label{lr:ray-profile-infinity-difference}
\normo{\wt V(t)-V(t)}_\infty
 \lsm_{p,B}\log t\,\normo{W_t-W}_\infty \lsm_{W,s} t^{d/4-s/2}\log t.
}
The standard difference estimate for $DF$ and~\eqref{lr:ray-profile-infinity-difference} give, for $h\in L^2$,
\EQn{\label{lr:linear-potential-difference}
\normo{(DF(\wt A(t))-DF(A_W(t)))[h]}_2
\lsm_{W,s}t^{-1-\nu}(\log t)^{\gamma_0}\normo{h}_2.
}

Rewrite~\eqref{lr:regularized-linear-equation} as
\EQ{
\brk{i\partial_t+\frac12\Delta}h
 =\lambda DF(A_W(t))[h]
   +f+\lambda\brkb{DF(\wt A(t))-DF(A_W(t))}[h].
}
For any fixed $\eta>0$, Corollary~\ref{lin:forced-estimate} and~\eqref{lr:linear-potential-difference} give, whenever $r\ge t\ge T_{W,s,\eta}$,
\EQn{\label{lr:original-flow-bound}
t^\eta\normo{h(t)}_2
 &\lsm_{W,s,\eta}
       r^\eta\normo{h(r)}_2
       +\int_t^r\tau^\eta\normo{f(\tau)}_2\,\dd \tau\\
 &\quad+\int_t^r
       \tau^{-1-\nu}(\log\tau)^{\gamma_0}
       \tau^\eta\normo{h(\tau)}_2\,\dd \tau.
}
Since $\tau^{-1-\nu}(\log\tau)^{\gamma_0}\in L^1([2,\infty))$, Gronwall's inequality gives~\eqref{lr:perturbed-flow-bound}.
\end{proof}

\begin{cor}[The linear final-state problem]\label{nl:inverse}
Let $W\in H^s(\R^d)$ with $d/2<s\le2$. Fix $\eta>0$ and take $T\ge T_{W,s,\eta}$ as in Proposition~\ref{lr:perturbed-flow}. Suppose that $f:[T,\infty)\to L^2(\R^d)$ is strongly measurable and
\EQ{
\int_T^\infty r^\eta\normo{f(r)}_{L^2}\,\dd r<\infty.
}
There exists a unique $v\in C([T,\infty);L^2)$ solving
\EQn{\label{eq:nl-forced-linear}
\brk{i\partial_t+\tfrac12\Delta}v
 =\lambda DF(\wt A)[v]+f
}
and satisfying
\EQn{\label{eq:nl-G}
\lim_{t\to\infty}t^\eta\normo{v(t)}_{L^2}=0.
}
Moreover, for every $R\ge T$,
\EQn{\label{eq:nl-inverse-S}
\normo{v}_{S([R,2R])}
 \lsm_{W,s,\eta} \int_R^\infty
       \brk{\frac rR}^\eta\normo{f(r)}_{L^2}\,\dd r.
}
\end{cor}

\begin{proof}
Let $\U(t,r)g$ denote the homogeneous solution given by Proposition~\ref{lr:perturbed-flow}, with terminal value $h(r)=g$. Set
\EQ{
v(t):=\int_t^\infty\U(t,r)[if(r)]\,\dd r.
}
By~\eqref{lr:perturbed-flow-bound}, this integral converges in $L^2$, and
\EQn{\label{eq:nl-inverse-L2}
t^\eta\normo{v(t)}_2
\lsm_{W,s,\eta}\int_t^\infty r^\eta\normo{f(r)}_2\,\dd r
\longrightarrow0\quad\text{as }t\to\infty.
}
For each $b>T$, the composition law gives
\EQ{
v(t)=\U(t,b)v(b)+\int_t^b\U(t,r)[if(r)]\,\dd r,
\qquad T\le t\le b.
}
The right-hand side is the unique solution on $[T,b]$ supplied by Proposition~\ref{lr:perturbed-flow} with terminal value $v(b)$. Hence $v\in C([T,\infty);L^2)$ solves~\eqref{eq:nl-forced-linear}. If $v,w$ are two solutions satisfying~\eqref{eq:nl-G}, estimate~\eqref{lr:perturbed-flow-bound} gives
\EQ{
\normo{v(t)-w(t)}_2
\lsm_{W,s,\eta}t^{-\eta}r^\eta\normo{v(r)-w(r)}_2
\longrightarrow0\quad\text{as }r\to\infty,
}
which proves uniqueness.

The Strichartz estimates, applied to~\eqref{eq:nl-forced-linear} with $t_0=R$, and~\eqref{lr:auxiliary-pointwise} give
\EQn{\label{eq:nl-free-S}
\normo{v}_{S([R,2R])}
&\lsm \normo{v(R)}_2
 +\normo{DF(\wt A)[v]}_{L^1([R,2R];L^2)}
 +\normo{f}_{L^1([R,2R];L^2)}\\
&\lsm_W \normo{v}_{L^\infty([R,2R];L^2)}
 +\normo{f}_{L^1([R,2R];L^2)}.
}
Together with~\eqref{eq:nl-inverse-L2}, this proves~\eqref{eq:nl-inverse-S}.
\end{proof}

\begin{lem}\label{nl:dyadic-forcing}
Let $W,s,\eta,T$ be as in Corollary~\ref{nl:inverse}. Fix $\mu>\eta$, $b\ge0$, and $D\ge0$. Suppose that $f:[T,\infty)\to L^2$ is strongly measurable and
\EQn{\label{eq:nl-dyadic-source}
\normo{f}_{L^1([2^jT,2^{j+1}T];L^2)}
 \le D (2^jT)^{-\mu}(\log(2^jT))^b,\qquad j\ge0.
}
Then equation~\eqref{eq:nl-forced-linear} has a unique solution $v\in C([T,\infty);L^2)$ satisfying~\eqref{eq:nl-G}, and for every $j\ge0$,
\EQn{\label{eq:nl-dyadic-target}
\normo{v}_{S([2^jT,2^{j+1}T])}
 \lsm_{W,s,\eta,\mu,b} D (2^jT)^{-\mu}(\log(2^jT))^b.
}
\end{lem}

\begin{proof}
Summing over the intervals $[2^{j+k}T,2^{j+k+1}T]$, we obtain
\EQ{
\int_{2^jT}^\infty\brk{\frac r{2^jT}}^\eta\normo{f(r)}_2\,\dd r
 \lsm_{\eta,\mu,b} D (2^jT)^{-\mu}(\log(2^jT))^b.
}
Taking $j=0$ gives $\int_T^\infty r^\eta\normo{f(r)}_2\,\dd r<\infty$, so Corollary~\ref{nl:inverse} gives the unique solution satisfying~\eqref{eq:nl-G}. Estimate~\eqref{eq:nl-inverse-S} now yields~\eqref{eq:nl-dyadic-target}.
\end{proof}

\section{Proof of the main theorem}\label{sec:nonlinear}

\subsection{Construction of the nonlinear correction}

\begin{prop}[Nonlinear final-state construction]\label{nl:construction}
Let $W\in H^s(\R^d)$ with $d/2<s\le2$. Let $A_W$ be the profile in~\eqref{eq:profile}. Then there exist $T_{W,s}\ge2$ and an $L^2$ solution $u$ of~\eqref{eq:nls} on $[T_{W,s},\infty)$, such that, for every $t\ge T_{W,s}$,
\EQn{\label{eq:nl-constructed-rate}
\normo{u-A_W}_{S([t,2t])}
 \lsm_{W,s} t^{-s/2}(\log t)^2.
}
No smallness condition is imposed on $W$.
\end{prop}

\begin{proof}
Let $\wt A$ be the profile in~\eqref{lr:profiles} and $\wt e$ its error in~\eqref{eq:error-exact}. Writing $u=\wt A+h$ in~\eqref{eq:nls} and using~\eqref{eq:profile-equation}, we obtain
\EQn{\label{eq:nl-perturbation}
\brk{i\partial_t+\frac12\Delta}h
=\lambda DF(\wt A(t))[h]+\lambda\mathcal N_{\wt A(t)}(h)-\wt e.
}
We solve this equation by the contraction mapping theorem. Fix
\EQ{
\eta:=\frac s4<\frac s2,\qquad T\ge T_{W,s,\eta},
}
where $T_{W,s,\eta}$ is given by Proposition~\ref{lr:perturbed-flow}. Define
\EQn{\label{eq:nl-X}
\normo{h}_{X_{s,T}}
 :=\sup_{j\ge0}\frac{(2^jT)^{s/2}}{(\log(2^jT))^2}
       \normo{h}_{S([2^jT,2^{j+1}T])}.
}
Here $X_{s,T}$ is the Banach space of functions in $C([T,\infty);L^2)\cap L^{\frac{2(d+2)}d}_{\mathrm{loc}}([T,\infty);L^{\frac{2(d+2)}d})$ with finite norm.

Proposition~\ref{prop:error} gives
\EQ{
\normo{\wt e}_{L^1([R,2R];L^2)}
 \lsm_{W,s} R^{-s/2}(\log R)^2.
}
Since $s/2>\eta$, dyadic summation gives $\int_T^\infty r^\eta\normo{\wt e(r)}_2\,\dd r<\infty$.

For $d=2,3$, H\"older's inequality yields, for $R\ge2$, every compact interval $I\subset[R,2R]$, and $h,k\in S(I)$,
\EQnn{\label{eq:nl-remainder-spacetime}
&\normo{\mathcal N_{\wt A}(h)-\mathcal N_{\wt A}(k)}_{L^1(I;L^2)}
 \notag\\
 &\quad\lsm_{p} R^{1/2}
 \brkb{\normo{h}_{L^{\frac{2(d+2)}d}(I\times\R^d)}^p
       +\normo{k}_{L^{\frac{2(d+2)}d}(I\times\R^d)}^p}
 \normo{h-k}_{L^{\frac{2(d+2)}d}(I\times\R^d)}.
}
Take $M\ge1$ and $h,k\in X_{s,T}$ with $\normo{h}_{X_{s,T}},\normo{k}_{X_{s,T}}\le M$. On each interval $[2^jT,2^{j+1}T]$, $j\ge0$, this gives
\EQnn{\label{eq:nl-remainder-dyadic}
 \normo{\mathcal N_{\wt A}(h)-\mathcal N_{\wt A}(k)}_{L^1([2^jT,2^{j+1}T];L^2)} \lsm_{p} M^p (2^jT)^{1/2-\frac{(d+2)s}{2d}}(\log(2^jT))^{2p+2}
                 \normo{h-k}_{X_{s,T}}.
}
The decay exponent satisfies
\EQ{
\mu:=\frac{(d+2)s}{2d}-\tfrac12>\tfrac s2>\eta,
}
because $s/d>1/2$. Taking $k=0$ in \eqref{eq:nl-remainder-dyadic} and summing as in \eqref{eq:nl-dyadic-source} therefore proves $\int_T^\infty r^\eta\normo{\mathcal N_{\wt A}(h)(r)}_2\,\dd r<\infty$.

For $d=1$, let $R\ge2$, $I\subset[R,2R]$ be a compact interval, and $h,k\in S(I)$. Since $\normo{\wt A}_{L^\infty(I;L^\infty)}\le BR^{-1/2}$, H\"older's inequality gives
\EQnn{\label{eq:nl-cubic-remainder-spacetime}
\normo{\mathcal N_{\wt A}(h)-\mathcal N_{\wt A}(k)}_{L^1(I;L^2)} 
&\lsm BR^{-1/2}
 \brkb{\normo{h}_{L^2(I;L^4)}+\normo{k}_{L^2(I;L^4)}}
 \normo{h-k}_{L^2(I;L^4)}
 \notag\\
&\quad+R^{1/2}
 \brkb{\normo{h}_{L^6(I;L^6)}^2+\normo{k}_{L^6(I;L^6)}^2}
 \normo{h-k}_{L^6(I;L^6)}
 \notag\\
&\lsm BR^{1/4}
 \brkb{\normo{h}_{S(I)}+\normo{k}_{S(I)}}\normo{h-k}_{S(I)}
 \notag\\
&\quad+R^{1/2}
 \brkb{\normo{h}_{S(I)}^2+\normo{k}_{S(I)}^2}\normo{h-k}_{S(I)},
}
where we used
\EQ{
\normo{f}_{L^2(I;L^4)}
 \le R^{3/8}\normo{f}_{L^\infty(I;L^2)}^{1/4}
             \normo{f}_{L^6(I;L^6)}^{3/4}
 \le R^{3/8}\normo{f}_{S(I)}.
}
Since $s>1/2$, we have $R^{1/4-s/2}(\log R)^2\lsm_s1$ for $R\ge2$. For $h,k$ in the same ball of $X_{s,T}$, the two terms in~\eqref{eq:nl-cubic-remainder-spacetime} therefore give, for every $j\ge0$,
\EQnn{\label{eq:nl-cubic-remainder-dyadic}
 \normo{\mathcal N_{\wt A}(h)-\mathcal N_{\wt A}(k)}_{L^1([2^jT,2^{j+1}T];L^2)}
\lsm_{W,s}M^2(2^jT)^{1/4-s}(\log(2^jT))^4
       \normo{h-k}_{X_{s,T}}.
}
The decay exponent satisfies $s-\tfrac14>\tfrac s2>\eta$. Taking $k=0$ and summing gives $\int_T^\infty r^\eta\normo{\mathcal N_{\wt A}(h)(r)}_2\,\dd r<\infty$ also in 1D.

With $\U(t,r)$ as in the proof of Corollary~\ref{nl:inverse}, define, for $h\in X_{s,T}$,
\EQn{\label{eq:nl-fixed-point}
(\mathcal Th)(t):=\int_t^\infty
 \U(t,r)\brkb{i\lambda\mathcal N_{\wt A}(h)(r)-i\wt e(r)}\,\dd r.
}
Corollary~\ref{nl:inverse} implies that $\mathcal Th$ is well defined. Then, we seek $h=\mathcal Th$ in $X_{s,T}$.

Since $\mathcal N_{\wt A}(0)=0$, Lemma~\ref{nl:dyadic-forcing}, applied to $-\wt e$ with $\mu=s/2$ and $b=2$, gives
\EQn{\label{eq:nl-error-X}
\normo{\mathcal T0}_{X_{s,T}}\le D_{W,s}.
}
For $\normo{h}_{X_{s,T}},\normo{k}_{X_{s,T}}\le M$, Corollary~\ref{nl:inverse} gives
\EQ{
 \normo{\mathcal Th-\mathcal Tk}_{S([2^jT,2^{j+1}T])} \lsm_{W,s}\int_{2^jT}^\infty
 \brk{\frac r{2^jT}}^\eta
 \normo{\mathcal N_{\wt A}(h)(r)-\mathcal N_{\wt A}(k)(r)}_2\,\dd r.
}
For $d=2,3$, Lemma~\ref{nl:dyadic-forcing} and~\eqref{eq:nl-remainder-dyadic}, with $\mu=\frac{(d+2)s}{2d}-1/2$ and $b=2p+2$, therefore yield
\EQn{\label{eq:nl-contraction-estimate}
\normo{\mathcal Th-\mathcal Tk}_{X_{s,T}}
 &\lsm_{W,s}M^p\sup_{j\ge0}
 (2^jT)^{1/2-s/d}(\log(2^jT))^{2p}
 \normo{h-k}_{X_{s,T}}\\
 &\lsm_{W,s}M^pT^{1/2-s/d}(\log T)^{2p}
 \normo{h-k}_{X_{s,T}}.
}
For $d=1$, applying the same lemma to~\eqref{eq:nl-cubic-remainder-dyadic} with $\mu=s-1/4$ and $b=4$ gives
\EQn{\label{eq:nl-cubic-contraction-estimate}
\normo{\mathcal Th-\mathcal Tk}_{X_{s,T}}
 &\lsm_{W,s}M^2\sup_{j\ge0}
 (2^jT)^{1/4-s/2}(\log(2^jT))^2
 \normo{h-k}_{X_{s,T}}\\
 &\lsm_{W,s}M^2T^{1/4-s/2}(\log T)^2
 \normo{h-k}_{X_{s,T}}.
}
Fix $M:=2D_{W,s}+1$. Since $s>d/2$, we may take $T$ sufficiently large that~\eqref{eq:nl-contraction-estimate} or~\eqref{eq:nl-cubic-contraction-estimate}, together with~\eqref{eq:nl-error-X}, gives
\EQ{
\normo{\mathcal Th}_{X_{s,T}}
 & \le\normo{\mathcal T0}_{X_{s,T}}
     +\normo{\mathcal Th-\mathcal T0}_{X_{s,T}} \le D_{W,s}+\tfrac12\normo{h}_{X_{s,T}}<M,\\
\normo{\mathcal Th-\mathcal Tk}_{X_{s,T}}
 &\le\tfrac12\normo{h-k}_{X_{s,T}}.
}
Thus $\mathcal T$ is a contraction on $\{h\in X_{s,T}:\normo{h}_{X_{s,T}}\le M\}$, and we obtain a unique fixed point $h=\mathcal Th$. By Corollary~\ref{nl:inverse} and~\eqref{eq:nl-fixed-point}, this fixed point solves~\eqref{eq:nl-perturbation}, so $u:=\wt A+h$ solves~\eqref{eq:nls}.

Finally, if $2^jT\le t<2^{j+1}T$, then 
\EQn{\label{lr:constructed-auxiliary-rate}
\normo{u-\wt A}_{S([t,2t])}
 \lsm_{W,s} t^{-s/2}(\log t)^2.
}
Proposition~\ref{lr:profile-comparison} yields
\EQ{
\normo{\wt A-A_W}_{S([t,2t])}
 \lsm_{W,s} t^{-s/2}\log t.
}
Combining these two estimates proves~\eqref{eq:nl-constructed-rate}. Since $u-A_W\in C([T,\infty);L^2)$, this also gives
\EQ{
\normo{u(t)-A_W(t)}_2\lsm_{W,s} t^{-s/2}(\log t)^2.
}
For the spacetime estimate,
\EQ{
\normo{u-A_W}_{L_{\tau,x}^{\frac{2(d+2)}d}([t,\infty)\times\R^d)}^{\frac{2(d+2)}d}
 &=\sum_{j=0}^\infty
   \normo{u-A_W}_{L_{\tau,x}^{\frac{2(d+2)}d}([2^jt,2^{j+1}t]\times\R^d)}^{\frac{2(d+2)}d}\\
 &\lsm_{W,s} t^{-\frac{s(d+2)}d}(\log t)^{\frac{4(d+2)}d}
   \sum_{j=0}^\infty 2^{-\frac{js(d+2)}d}(1+j)^{\frac{4(d+2)}d}\\
 &\lsm_{W,s} t^{-\frac{s(d+2)}d}(\log t)^{\frac{4(d+2)}d}.
}
Combining it with the pointwise $L^2$ bound proves~\eqref{eq:main-rate}.
\end{proof}

\subsection{Uniqueness in the \texorpdfstring{$L^2$}{L2} class and the wave operator}\label{sec:uniqueness}

To prove uniqueness, we first recover the spacetime decay from the assumed $L^2$ decay.

\begin{lem}[Spacetime estimate]\label{nl:bootstrap}
Let $W\in H^s(\R^d)$ with $d/2<s\le2$. Fix $d/4<\beta<s/2$ and $T_0>0$. Suppose that $u$ is an $L^2$ solution of~\eqref{eq:nls} on $[T_0,\infty)$ and that, for some finite $M\ge0$,
\EQn{\label{eq:nl-assumed-L2-rate}
\normo{u(t)-A_W(t)}_{L^2}\le M t^{-\beta}
}
for all sufficiently large $t$. Then there exists a sufficiently large $T_1\ge T_0$ such that
\EQn{\label{eq:nl-automatic-S-rate}
\normo{u-A_W}_{S([R,2R])}\lsm_{W,s,M,\beta}R^{-\beta},\qquad R\ge T_1.
}
\end{lem}

\begin{proof}
Set $h:=u-\wt A$. Since $\beta<s/2$, Proposition~\ref{lr:profile-comparison} gives, for all sufficiently large $t$,
\EQn{\label{lr:bootstrap-auxiliary-L2}
\normo{h(t)}_2\le(M+1)t^{-\beta}.
}
For $R<r\le2R$, set $Y_R(r):=\normo{h}_{S([R,r])}$. Applying the Strichartz estimates to~\eqref{eq:nl-perturbation}, and using~\eqref{lr:auxiliary-pointwise}, \eqref{lr:bootstrap-auxiliary-L2}, and Proposition~\ref{prop:error}, we obtain
\EQ{
Y_R(r)\le C_0R^{-\beta}
 +C_d\normo{\mathcal N_{\wt A}(h)}_{L^1([R,r];L^2)},
}
where $C_0=C_0(W,s,M,\beta)\ge M+1$. Estimates~\eqref{eq:nl-remainder-spacetime} and~\eqref{eq:nl-cubic-remainder-spacetime}, with $k=0$, give
\EQn{\label{eq:nl-bootstrap-inequality}
Y_R(r)\le C_0R^{-\beta}+C_pR^{1/2}Y_R(r)^{p+1},
\qquad d=2,3,
}
and
\EQn{\label{eq:nl-cubic-bootstrap-inequality}
Y_R(r)\le C_0R^{-\beta}
 +CBR^{1/4}Y_R(r)^2+CR^{1/2}Y_R(r)^3,
\qquad d=1.
}
Since $p\beta>1/2$ for $d=2,3$ and $\beta>1/4$ for $d=1$, these estimates imply that, for all sufficiently large $R$,
\EQ{
Y_R(r)\le2C_0R^{-\beta}
\quad\Longrightarrow\quad
Y_R(r)\le\tfrac32C_0R^{-\beta}.
}
The function $Y_R$ extends continuously to $r=R$, with $Y_R(R)=\normo{h(R)}_2\le C_0R^{-\beta}$. Hence the bootstrap argument gives $Y_R(2R)\le2C_0R^{-\beta}$. Proposition~\ref{lr:profile-comparison} then yields~\eqref{eq:nl-automatic-S-rate}.
\end{proof}

\begin{prop}[Final-state uniqueness]\label{nl:uniqueness}
Let $W\in H^s(\R^d)$ with $d/2<s\le2$. Let $T_0>0$, and let $u_1,u_2$ be $L^2$ solutions of~\eqref{eq:nls} on $[T_0,\infty)$. Suppose that, for each $j=1,2$, there exists $\beta_j>d/4$ such that
\EQ{
\normo{u_j(t)-A_W(t)}_{L^2}=O(t^{-\beta_j})
 \quad\text{as }t\longrightarrow\infty.
}
Then $u_1=u_2$ on their common interval of existence.
\end{prop}

\begin{proof}
Let $h_j:=u_j-\wt A$ and $k:=h_1-h_2=u_1-u_2$. Equation~\eqref{eq:nl-perturbation} gives
\EQ{
 \brk{i\partial_t+\tfrac12\Delta}k
 =\lambda DF(\wt A)[k]
  +\lambda\brkb{\mathcal N_{\wt A}(h_1)
                 -\mathcal N_{\wt A}(h_2)}.
}
Choose
\EQ{
\frac d4<\beta<\min\{s/2,\beta_1,\beta_2\},
 \qquad 0<\eta<\beta,
}
By Lemma~\ref{nl:bootstrap} and Proposition~\ref{lr:profile-comparison}, after increasing $T$ we have
\EQn{\label{eq:nl-Xbeta}
\normo{h_j}_{X_\beta(T)}\le M_j,
 \qquad
 \normo{v}_{X_\beta(T)}
 :=\sup_{\ell\ge0}(2^\ell T)^\beta
       \normo{v}_{S([2^\ell T,2^{\ell+1}T])}.
}
Take $T\ge T_{W,s,\eta}$. Since $\eta<\beta$, \eqref{eq:nl-Xbeta} gives
\EQ{
t^\eta\normo{k(t)}_2
 \lsm_\beta(M_1+M_2)t^{\eta-\beta}\longrightarrow0
 \quad\text{as }t\to\infty.
}

In $d=2,3$, for $R=2^\ell T$, \eqref{eq:nl-remainder-spacetime} and \eqref{eq:nl-Xbeta} yield
\EQnn{\label{eq:nl-difference-source}
 \normo{\mathcal N_{\wt A}(h_1)-\mathcal N_{\wt A}(h_2)}_{L^1([R,2R];L^2)} \lsm_{p} (M_1^p+M_2^p)\normo{k}_{X_\beta(T)}
           R^{1/2-(p+1)\beta}.
}
Here $
(p+1)\beta-\tfrac12>\beta>\eta$,
because $p\beta>1/2$.

In 1D case, $R^{1/2-3\beta}\le R^{1/4-2\beta}$ for $R\ge2$, so~\eqref{eq:nl-cubic-remainder-spacetime} gives
\EQnn{\label{eq:nl-cubic-difference-source}
 \normo{\mathcal N_{\wt A}(h_1)-\mathcal N_{\wt A}(h_2)}_{L^1([R,2R];L^2)} \lsm_{W,M_1,M_2}R^{1/4-2\beta}\normo{k}_{X_\beta(T)}.
}
The decay exponent satisfies 
$2\beta-\tfrac14>\beta>\eta$, because $\beta>1/4$. Summing~\eqref{eq:nl-difference-source} or~\eqref{eq:nl-cubic-difference-source} over dyadic intervals proves
\EQ{
\int_T^\infty r^\eta
 \normo{\mathcal N_{\wt A}(h_1)(r)-\mathcal N_{\wt A}(h_2)(r)}_2\,\dd r<\infty.
}
Corollary~\ref{nl:inverse} therefore identifies $k$ with the solution of the linear final-state problem with forcing $\lambda(\mathcal N_{\wt A}(h_1)-\mathcal N_{\wt A}(h_2))$.

We now estimate $k$ on $[T',\infty)$, where $T':=2^mT$. The bounds in~\eqref{eq:nl-Xbeta} hold there with the same constants $M_j$. For $d=2,3$, Lemma~\ref{nl:dyadic-forcing}, with $\mu=(p+1)\beta-1/2$ and $b=0$, gives
\EQ{
\normo{k}_{X_\beta(T')}
 \lsm_{W,s,\beta,\eta,M_1,M_2}(T')^{1/2-p\beta}
       \normo{k}_{X_\beta(T')}.
}
When $d=1$, applying Lemma~\ref{nl:dyadic-forcing} to~\eqref{eq:nl-cubic-difference-source} with $\mu=2\beta-1/4$ and $b=0$ gives
\EQ{
\normo{k}_{X_\beta(T')}
 \lsm_{W,s,\beta,\eta,M_1,M_2}(T')^{1/4-\beta}
       \normo{k}_{X_\beta(T')}.
}
The norm on the right is finite by \eqref{eq:nl-Xbeta}. For $d=2,3$, we have $1/2-p\beta<0$. For $d=1$, we have $1/4-\beta<0$. Choosing $m$ large enough gives $k=0$ on $[T',\infty)$. Uniqueness for the $L^2$ Cauchy problem~\cite{Tsutsumi1987} extends this equality to the entire common interval.
\end{proof}

\begin{cor}\label{nl:global}
Let $u_+\in \F H^s(\R^d)$ with $d/2<s\le2$, and set $W:=\wh{u_+}$. Then the solution $u$ constructed in Proposition~\ref{nl:construction} extends uniquely to a global $L^2$ solution.
The modified wave operator
\EQ{
\Om_+:\F H^s(\R^d)\longrightarrow L^2(\R^d),
 \qquad \Om_+(u_+):=u(0),
}
is well defined in the asymptotic class of Theorem~\ref{thm:main}.
\end{cor}

\begin{proof}
By Tsutsumi's result~\cite{Tsutsumi1987}, the datum $u(T_{W,s})$ generates a unique global $L^2$ solution. Thus the uniqueness implies that $\Om_+(u_+)=u(0)$ is well defined on $\F H^s$, which completes the proof of Theorem~\ref{thm:main}.
\end{proof}

\end{document}